\documentclass[11pt, a4paper]{article}
\usepackage{times}
\usepackage{a4wide}
\usepackage[dvips, hyperindex]{hyperref}
\usepackage[british]{babel}
\usepackage{enumerate, longtable}
\usepackage{amsmath, amscd, amsfonts, amssymb, latexsym, theorem, comment}
\usepackage{xypic}
\usepackage[T1]{fontenc}
\usepackage[latin1]{inputenc}

{\theorembodyfont{\itshape}   \newtheorem{thm}{Theorem}[section]}
{\theorembodyfont{\itshape}   \newtheorem{lem}[thm]{Lemma}}
{\theorembodyfont{\itshape}   \newtheorem{defi}[thm]{Definition}}
{\theorembodyfont{\itshape}   }
{\theorembodyfont{\itshape}   }
{\theorembodyfont{\itshape}   \newtheorem{prop}[thm]{Proposition}}
{\theorembodyfont{\itshape}   \newtheorem{cor}[thm]{Corollary}}
{\theorembodyfont{\itshape}   }
{\theorembodyfont{\itshape}   }
{\theorembodyfont{\itshape}   }
{\theorembodyfont{\itshape}   \newtheorem{question}[thm]{Question}}
{\theorembodyfont{\itshape}   \newtheorem{notation}[thm]{Notation}}

\newcommand{\CC}{\mathbb{C}}
\newcommand{\QQ}{\mathbb{Q}}
\newcommand{\FF}{\mathbb{F}}
\newcommand{\ZZ}{\mathbb{Z}}
\newcommand{\PP}{\mathbb{P}}

\newcommand{\TT}{\mathbb{T}}

\newcommand{\NN}{\mathbb{N}}

\newcommand{\Qbar}{{\overline{\QQ}}}

\newcommand{\Fpbar}{{\overline{\FF}_p}}
\newcommand{\chibar}{{\overline{\chi}}}

\newcommand{\Hom}{{\rm Hom}}

\newcommand{\GL}{\mathrm{GL}}
\newcommand{\SL}{\mathrm{SL}}

\DeclareMathOperator{\Gal}{Gal}

\DeclareMathOperator{\Frob}{Frob}
\DeclareMathOperator{\Tr}{Tr}
\DeclareMathOperator{\Det}{Det}

\newcommand{\Ind}{{\rm Ind}}

\newcommand{\pf}{{\bf Proof. }}
\newcommand{\qed}{\hspace* {.5cm} \hfill $\Box$}

\newcommand{\mat}[4]{
 \left(  \begin{smallmatrix} #1 & #2 \\ #3 & #4 \end{smallmatrix} \right)}

\newcommand{\cO}{\mathcal{O}}
\newcommand{\Qpbar}{{\overline{\QQ}_p}}
\newcommand{\Ftbar}{{\overline{\FF}_2}}

\newcommand{\fP}{\mathfrak{P}}
\newcommand{\fQ}{\mathfrak{Q}}
\newcommand{\fR}{\mathfrak{R}}
\newcommand{\cM}{\mathcal{M}}
\newcommand{\cCM}{\mathcal{CM}}
\newcommand{\cB}{\mathcal{B}}
\newcommand{\fm}{\mathfrak{m}}
\newcommand{\fp}{\mathfrak{p}}
\newcommand{\rk}{{\rm rk}}
\newcommand{\Eis}{{\rm Eis}}

\newcommand{\magma}[1]{{\bf\textsl{#1}}}
\newcommand{\magmaout}[1]{{\small\bf\textit{#1}}}
\newcommand{\intr}[1]{\medskip{\noindent \magma{#1} \\}}

\newcommand{\m}{\mathfrak{m}}
\newcommand{\hecke}[1]{\TT_\ZZ(S_{k}(#1,\chi))}

\newcommand{\heckeoneNk}[2]{\TT_\ZZ(S_{#2}(\Gamma_1(#1)))}

\newcommand{\heckeNktwo}[1]{\TT_\ZZ(S_2(#1,\chi_{\rm triv}))}

\title{On the failure of the Gorenstein property\\ 
for Hecke algebras of prime weight}
\author{L. J. P. Kilford and Gabor Wiese}

\begin{document}

\maketitle

\begin{abstract}
  In this article we report on extensive calculations concerning the
  Gorenstein defect for Hecke algebras of spaces of modular forms of
  prime weight~$p$ at maximal ideals of residue characteristic~$p$
  such that the attached mod~$p$ Galois representation is unramified
  at~$p$ and the Frobenius at~$p$ acts by scalars. The results lead us
  to the ask the question whether the Gorenstein defect and the
  multplicity of the attached Galois representation are always equal
  to~$2$. We review the literature on the failure of the Gorenstein
  property and multiplicity one, discuss in some detail a very
  important practical improvement of the modular symbols algorithm
  over finite fields and include precise statements on the
  relationship between the Gorenstein defect and the multiplicity of
  Galois representations.

\end{abstract}

\tableofcontents

\section{Introduction}
\label{section-introduction}

In Wiles' proof of Fermat's Last Theorem (see \cite{wiles}) an
essential step was to show that certain Hecke algebras are Gorenstein
rings. Moreover, the Gorenstein property of Hecke algebras is
equivalent to the fact that Galois representations appear on certain
Jacobians of modular curves precisely with multiplicity one. This
article is concerned with the Gorenstein property and with the
multiplicity one question. We report previous work and exhibit many
new examples where multiplicity one and the Gorenstein property fail.
We compute the multiplicty in these cases.
Moreover, we ask the question suggested by our computations whether it
is always equal to two if it fails.  \medskip

We first have to introduce some notation. For integers $N \ge 1$
and $k \ge 2$ and a Dirichlet character $\chi: (\ZZ/N\ZZ)^\times \to \CC^\times$
we let $S_k(\Gamma_1(N))$ be the $\CC$-vector space of holomorphic
cusp forms on $\Gamma_1(N)$ of weight~$k$ and $S_k(N,\chi)$ the
subspace on which the diamond operators act through the
character~$\chi$. We now introduce some extra notation for Hecke algebras over 
specified rings.
\begin{notation}[Notation for Hecke algebras]
Whenever $S \subseteq R$ are rings and
$M$ is an $R$-module on which the Hecke and diamond operators act,
we let $\TT_S(M)$ be the $S$-subalgebra inside the $R$-endomorphism
ring of~$M$ generated by the Hecke and the diamond operators.
If $\phi: S \to S'$ is a ring homomorphism, we let
$\TT_\phi (M) := \TT_S(M) \otimes_S S'$ or with $\phi$ understood
$\TT_{S \to S'} (M)$.
\end{notation}
We will mostly be dealing with the Hecke algebras
$\TT_\ZZ(S_k(\Gamma_1(N)))$ and $\TT_{\ZZ[\chi]}(S_k(N,\chi))$,
their completions $\TT_{\ZZ \to \ZZ_p}(S_k(\Gamma_1(N)))$ and
$\TT_{\cO\to \widehat{\cO}}(S_k(N,\chi))$, as well as their reductions
$\TT_{\ZZ \to \FF_p}(S_k(\Gamma_1(N)))$ and $\TT_{\cO \to
\FF}(S_k(N,\chi))$.  Here, $p$ is a prime and $\cO = \ZZ[\chi]$ is the
smallest subring of~$\CC$ containing all values of~$\chi$,
$\widehat{\cO}$ is a completion and $\cO \twoheadrightarrow \FF$ is a
reduction modulo a prime above~$p$. In Section~\ref{modular-symbols}
the reduced Hecke algebras are identified with Hecke algebras of
mod~$p$ modular forms, which are closely related to Hecke algebras of
Katz modular forms over finite fields (see
Section~\ref{section-relation}).

We choose a holomorphic cuspidal Hecke eigenform as the starting point
of our discussion and treatment. So let $f \in S_k(N,\chi) \subseteq
S_k(\Gamma_1(N))$ be an eigenform for all Hecke and diamond operators.
It (more precisely, its Galois conjugacy class) corresponds to minimal
ideals, both denoted by $\fp_f$, in each of the two Hecke algebras
$\TT_\ZZ(S_k(\Gamma_1(N)))$ and $\TT_{\ZZ[\chi]}(S_k(N,\chi))$.  We
also choose maximal ideals $\m=\fm_f$ containing~$\fp_f$ of residue
characteristic~$p$ again in each of the two. Note that the residue
fields are the same in both cases.

By work of Shimura and Deligne, one can associate to~$f$ (more
precisely, to~$\m$) a continuous odd semi-simple Galois representation
\[
\rho_f = \rho_{\m_f} = \rho_\m: \Gal(\Qbar/\QQ) \rightarrow \GL_2(\hecke{N}/\m)
\]
unramified outside $Np$ satisfying~$\Tr(\rho_\m(\Frob_l)) \equiv T_l
\mod \m$ and $\Det(\rho_m(\Frob_l)) \equiv l^{k-1}\chi(l) \mod \m$ for all
primes~$l \nmid Np$. In the case of weight~$k=2$ and level~$N$, the
representation~$\rho_m$ can be naturally realised on the $p$-torsion
points of the Jacobian of the modular curve $X_1(N)_\QQ$. The algebra
$\TT_{\ZZ \to \FF_p}(S_2(\Gamma_1(N)))$ acts naturally on
$J_1(N)_\QQ(\Qbar)[p]$ and we can form the Galois module $J_1(N)_\QQ
(\Qbar)[\m] = J_1(N)_\QQ(\Qbar)[p][\widetilde{\m}]$ with
$\widetilde{\m}$ the maximal ideal of $\TT_{\ZZ \to
\FF_p}(S_2(\Gamma_1(N)))$ which is the image of~$\m$ under the natural
projection. Supposing that $\rho_m$ is absolutely irreducible, the
main result of~\cite{blr} shows that the Galois representation
$J_1(Np)_\QQ(\Qbar)[\m]$ is isomorphic to a direct sum of $r$~copies
of~$\rho_m$ for some integer $r \ge 1$, which one calls the {\em
multiplicity of~$\rho_m$} (cf.~\cite{ribet-stein}). One says that
$\rho_m$ is a {\em multiplicity one representation} or {\em satisfies
multiplicity one}, if $r=1$. See \cite{mazur} for a similar definition
for $J_0(N)$ and Prop.~\ref{nulleins} for a comparison.

The notion of multiplicity can be naturally extended to Galois
representations attached to eigenforms~$f$ of weights $3 \le k \le p+1$
for $p \nmid N$. This is accomplished by a result of Serre's which
implies the existence of a maximal ideal $\m_2 \subset
\heckeoneNk{Np}{2}$ such that $\rho_{\m_f} \cong \rho_{\m_2}$
(see Prop.~\ref{totwo}). One hence obtains the notion of multiplicity
(on $J_1(Np)$) for the representation~$\rho_{m_f}$ by defining it as
the multiplicity of~$\rho_{\m_2}$.  Moreover, when allowing twists by
the cyclotomic character, it is even possible to treat arbitrary
weights.  The following theorem summarises results on when the
multiplicity in the above sense is known to be one.

\begin{thm}[Mazur~\cite{mazur}, Edixhoven~\cite{edixhoven-serre},
  Gross~\cite{gross}, Buzzard~\cite{buzzard_appendix}]
\label{gorenstein-theorem}
Let~$\rho_\m$ be a representation associated with a modular cuspidal
eigenform~$f \in S_{k}(N,\chi)$ and let~$p$ be the residue
characteristic of~$\m$. Suppose that $\rho_\m$ is absolutely
irreducible and that $p$ does not divide~$N$. If either
\begin{enumerate}
\itemsep=0cm plus 0pt minus 0pt
\item $2 \le k \le p-1$, or

\item $k=p$ and $\rho_\m$ is ramified at~$p$, or

\item $k=p$ and $\rho_\m$ is unramified at~$p$ and~$\rho_\m(\Frob_p)$ is not scalar,

\end{enumerate}
then the multiplicity of~$\rho_\m$ is one.
\end{thm}

This theorem is composed of Lemma~15.1 from Mazur~\cite{mazur},
Theorem~9.2 from Edixhoven~\cite{edixhoven-serre}, Proposition~{12.10}
from Gross~\cite{gross} and Theorem~6.1 from
Buzzard~\cite{buzzard_appendix}.
\medskip

The following theorem by the second author (\cite{nongor}, Corollary~4.4)
tells us when the multiplicity is not one.

\begin{thm}\label{nongor}
  Let $\rho_\m$ as in the previous theorem. Suppose $k=p$ and that
  $\rho_\m$ is unramified at~$p$. If $p=2$, assume also that a Katz
  cusp form over~$\FF_2$ of weight~$1$ on~$\Gamma_1(N)$ exists which
  gives rise to~$\rho_\m$.

  If $\rho_\m(\Frob_p)$ is a scalar matrix, then the multiplicity
  of~$\rho_\m$ is bigger than~$1$.
\end{thm}

In Section~\ref{section-relation} we explain how the Galois
representation $J_1(Np)_\QQ(\Qbar)[\m]$ is related to the different
Hecke algebras evoked above and see in many cases of interest a
precise relationship between the geometrically defined term of {\em
multiplicity} and the {\em Gorenstein defect} of these algebras. The
latter can be computed explicitly, which is the subject of the present
article. We now give the relevant definitions.

\begin{defi}[The Gorenstein property]
  Let~$A$ be a local Noetherian ring with maximal ideal~$\m$.
  Suppose first that the Krull dimension of~$A$ is zero, i.e.\ that
  $A$ is Artinian. We then define the \emph{Gorenstein defect} of~$A$
  to be the minimum number of $A$-module generators of the annihilator
  of~$\m$ (i.e.\ $A[\m]$) minus~1; equivalently, this is
  the~$A/\m$-dimension of the annihilator of~$\m$ minus~1. We say that
  $A$ is \emph{Gorenstein} if its Gorenstein defect is~0, and
  \emph{non-Gorenstein} otherwise.
  If the Krull dimension of~$A$ is positive, we inductively call~$A$
  {\em Gorenstein}, if there exists a non-zero-divisor $x \in \fm$
  such that $A/(x)$ is a Gorenstein ring of smaller Krull dimension
  (see~\cite{Eisenbud}, p.~532; note that our definition implies that
  $A$ is Cohen-Macaulay).
  A (not necessarily local) Noetherian ring is said to be
  \emph{Gorenstein} if and only if all of its localisations at its
  maximal ideals are Gorenstein.
\end{defi}

We will for example be interested in the Gorenstein property of
$\TT_\ZZ(S_k(\Gamma_1(N)))_\m$. Choosing $x=p$ in the definition, we
see that this is equivalent to the Gorenstein defect of the finite
dimensional $\FF_p$-algebra $\TT_{\ZZ\to\FF_p}(S_k(\Gamma_1(N)))_\m$
being zero. Whenever we refer to the Gorenstein defect of the former
algebra (over~$\ZZ$), we mean the one of the latter. Our computations
will concern the Gorenstein defect of
$\TT_{\cO\to\FF}(S_k(\Gamma_1(N),\chi))_\m$. See
Section~\ref{section-relation} for a comparison with the one
not involving a character. It is important to remark that the
Gorenstein defect of a local Artin algebra over a field does not
change after passing to a field extension and taking any of the
conjugate local factors.

We illustrate the definition by an example. The algebra
$k[x,y,z]/(x^2,y^2,z^2,xy,xz,yz)$ for a field~$k$ is Artinian and
local with maximal ideal $\m:=(x,y,z)$ and the annihilator of~$\m$
is~$\m$ itself, so the Gorenstein defect is~$3-1=2$. We note that this
particular case does occur in nature; a
localisation~$\TT_{\ZZ\to\FF_2}(S_2(\Gamma_0(431)))_\m$ at one maximal
ideal is isomorphic to this, with~$k=\FF_2$ (see~\cite{emerton}, the
discussion just before Lemma~6.6). This example can also be verified
with the algorithm presented in this paper.
\medskip

We now state a translation of Theorem~\ref{gorenstein-theorem}
in terms of Gorenstein defects, which is immediate from the
propositions in Section~\ref{section-relation}.

\begin{thm}\label{thgor}
Assume the set-up of Theorem~\ref{gorenstein-theorem} and that
one of 1.,\ 2.,\ or 3.\ is satisfied. We use notations as in the discussion
of multiplicities above.

If $k=2$, then $\TT_\ZZ(S_2(\Gamma_1(N)))_\m$ is a Gorenstein ring.

If $k\ge3$, then $\TT_\ZZ(S_2(\Gamma_1(Np)))_{\m_2}$ is, too.
Supposing in addition that $\m$ is ordinary (i.e.\ $T_p \not\in \m$), then
also $\TT_\ZZ(S_k(\Gamma_1(N)))_\m$ is Gorenstein. If, moreover, $p \ge 5$
or $\rho_\m$ is not induced from $\QQ(\sqrt{-1})$ (if $p=2$)
or $\QQ(\sqrt{-3})$ (if $p=3$), then
$\TT_{\cO\to\FF}(S_k(N,\chi))_\m$ is Gorenstein as well.
\qed
\end{thm}

We now turn our attention to computing the Gorenstein defect and the
multiplicity in the case when it is known not to be one.

\begin{cor}\label{corfailure}
Let~$\rho_\m$ be a representation associated with a cuspidal
eigenform~$f \in S_{p}(N,\chi)$ with~$p$ the residue
characteristic of~$\m$. Assume that $\rho_\m$ is absolutely irreducible, 
unramified at~$p$ such that $\rho_\m(\Frob_p)$ is a scalar matrix.
Let $r$ be the multiplicity of~$\rho_\m$ and
$d$ the Gorenstein defect of any of
$\TT_{\cO\to\FF}(S_k(N,\chi))_\m$,
$\TT_{\ZZ\to\FF_p}(S_k(\Gamma_1(N)))_\m$ or
$\TT_{\ZZ\to\FF_p}(S_2(\Gamma_1(Np)))_{\m_2}$.

Then the relation $d = 2r-2$ holds.
\end{cor}

\pf
The equality of the Gorenstein defects and
the relation with the multiplicity are proved in
Section~\ref{section-relation}, noting
that $\m$ is ordinary, as
$a_p(f)^2 = \chi(p) \neq 0$ by \cite{gross}, p.~487.
\qed

\subsection{Previous results on the failure of multiplicity one or the Gorenstein property}
\label{previous-results}

Prior to the present work and the article~\cite{nongor},
there have been some investigations into when Hecke algebras fail to
be Gorenstein. In~\cite{kilford-nongorenstein}, the first author
showed, using {\sc Magma}~\cite{magma}, that the Hecke
algebras~$\heckeNktwo{431}$, $\heckeNktwo{503}$
and~$\heckeNktwo{2089}$ are not Gorenstein by explicit computation of
the localisation of the Hecke algebra at a suitable maximal ideal
above~2, and in~\cite{ribet-stein}, it is shown
that~$\heckeNktwo{2071}$ is not Gorenstein in a similar fashion. These
examples were discovered by considering elliptic curves~$E/\QQ$ such that
in the ring of integers of~$\QQ(E[2])$ the prime ideal~$(2)$ splits
completely, and then doing computations with {\sc Magma}.

There are also some results in the literature on the failure of
multiplicity one within the torsion of certain Jacobians.
In~\cite{agashe-ribet-stein}, the following theorem is proved:
\begin{thm}[Agashe-Ribet-Stein~\cite{agashe-ribet-stein}, Proposition~5.1]
  Suppose that~$E$ is an optimal elliptic curve over~$\QQ$ of
  conductor~$N$, with congruence number~$r_E$ and modular degree~$m_E$
  and that~$p$ is a prime such that~$p | r_E$ but~$p \nmid m_E$.
  Let~$\mathfrak{m}$ be the annihilator in~$\heckeNktwo{N}$ of~$E[p]$.
  Then multiplicity one fails for~$\mathfrak{m}$.
\end{thm}
They give a table of examples; for instance, $\heckeNktwo{54}$ does
not satisfy multiplicity one at some maximal ideal above~$3$. It is
not clear whether this phenomenon occurs infinitely often.

In~\cite{ribet-festschrift}, it is shown that the mod~$p$ multiplicity
of a certain representation in the Jacobian of the Shimura curve
derived from the rational quaternion algebra of 
discriminant~$11 \cdot 193$ is~2; this result inspired the calculations
in~\cite{kilford-nongorenstein}.

Let us finally mention that for $p=2$ precisely the Galois
representations~$\rho$ with image equal to the dihedral group~$D_3$
come from an elliptic curve over~$\QQ$. For, on the one hand, we only
need observe that $D_3=\GL_2(\FF_2)$. On the other hand, any
$S_3$-extension~$K$ of the rationals can be obtained as the splitting
field of an irreducible integral polynomial $f=X^3+aX+b$. The
$2$-torsion of the elliptic curve $E: Y^2=f$ consists precisely of the
three roots of~$f$ and the point at infinity. So, the field generated
over $\QQ$ by the $2$-torsion of~$E$ is~$K$.

\subsection{New results}

Using the modular symbols algorithm over finite fields with an
improved stop criterion (see Section~\ref{modular-symbols}), we
performed computations in {\sc Magma} concerning the Gorenstein
defect of Hecke algebras of cuspidal modular forms of prime
weights~$p$ at maximal ideals of residue characteristic~$p$ in the case
of Theorem~\ref{nongor}. All of our 384
examples have Gorenstein defect equal to~$2$ and hence their
multiplicity is~$2$.

We formulate part of our computational findings as a theorem.

\begin{thm}
  For every prime $p < 100$ there exists a prime $N \neq p$ and a
  Dirichlet character~$\chi$ such that the Hecke algebra
  $\TT_{\ZZ[\chi]\to\FF}(S_p(N,\chi))$ has Gorenstein defect~$2$ at
  some maximal ideal~$\m$ of residue characteristic~$p$. The
  corresponding Galois representation~$\rho_\m$ appears with
  multiplicity two on the $p$-torsion of the Jacobian $J_1(Np)_\QQ(\Qbar)$
  if $p$ is odd, respectively $J_1(N)_\QQ(\Qbar)$ if~$p=2$.
\end{thm}

Our computational results are discussed in more detail in
Section~\ref{computational-results}.

\subsection{A question}

\begin{question}\label{ourquestion}
Let $p$ be a prime. Let $f$ be a normalised cuspidal modular
eigenform of weight~$p$, prime level $N \neq p$ for some Dirichlet
character $\chi$. Let $\rho_f : G_\QQ \to \GL_2(\Fpbar)$ be the
modular Galois representation attached to~$f$. We assume that
$\rho_f$ is irreducible and unramified at~$p$ and that
$\rho_f(\Frob_p)$ is a scalar matrix.

Write $\TT_{\FF}$ for $\TT_{\ZZ[\chi]\to\FF}(S_p(N,\chi))$. Recall
that this notation stands for the tensor product over $\ZZ[\chi]$ of a
residue field $\FF/\FF_p$ of $\ZZ[\chi]$ by the $\ZZ[\chi]$-algebra
generated inside the endomorphism algebra of~$S_p(N,\chi)$ by the
Hecke operators and by the diamond operators.  Let $\mathfrak{m}$ be
the maximal ideal of $\TT_{\FF}$ corresponding to~$f$.

Is the Gorenstein defect of the Hecke algebra~$\TT_{\FF}$ 
localised at ${\mathfrak{m}}$,
denoted by $\TT_{\mathfrak{m}}$, always equal to~$2$?

Equivalently, is the multiplicity of the Galois representation
attached to~$f$ always equal to~$2$?
\end{question}

This question was also raised both by Kevin Buzzard and James Parson in
communications to the authors.

\section{Relation between multiplicity and Gorenstein defect}
\label{section-relation}

In this section we collect results, some of which are well-known, on
the multiplicity of Galois representations, the Gorenstein defect and
relations between the two. Whereas the mod~$p$ modular symbols
algorithm naturally computes mod~$p$ modular forms (see
Section~\ref{modular-symbols}), this rather geometrical section uses
(mostly in the references) the theory of Katz modular forms over
finite fields (see e.g.~\cite{edixhoven-boston}). If $N \ge 5$ and $k\ge 2$,
the Hecke algebra $\TT_{\ZZ \to \FF_p}(S_k(\Gamma_1(N)))$ is both the
Hecke algebra of mod~$p$ cusp forms of weight~$k$ on~$\Gamma_1(N)$ and
the Hecke algebra of the corresponding Katz cusp forms over~$\FF_p$.
However, in the presence of a Dirichlet character 
$\TT_{\ZZ[\chi] \to \FF}(S_k(N,\chi))$ only has an interpretation as
the Hecke algebra of the corresponding mod~$p$ cusp forms and there
may be differences with the respective Hecke algebra for Katz forms
(see Carayol's Lemma, Prop.~1.10 of~~\cite{edixhoven-boston}).

We start with the well-known result in weight~$2$ (see e.g.\
\cite{mazur}, Lemma 15.1) that multiplicity one implies that the
corresponding local Hecke factor is a Gorenstein ring.

\begin{prop}
  Let $\fm$ be a maximal ideal of~$\TT := \heckeoneNk{N}{2}$ of residue
  characteristic~$p$ which may divide~$N$. Denote by $\widetilde{\m}$
  the image of $\m$ in 
  $\TT_{\FF_p} := \TT \otimes_\ZZ \FF_p = \TT_{\ZZ\to\FF_p}(S_2(\Gamma_1(N)))$. 
  Suppose that the Galois representation~$\rho_\m$ is irreducible and
  satisfies multiplicity one (see Section~\ref{section-introduction}).

  Then as $\TT_{\FF_p,\widetilde{\m}}$-modules one has 
  $$ J_1(N)_\QQ(\Qbar)[p]_{\widetilde{\m}} \cong 
     \TT_{\FF_p,\widetilde{\m}} \oplus \TT_{\FF_p,\widetilde{\m}}$$ 
  and the localisations $\TT_\m$ and $\TT_{\FF_p,\widetilde{\m}}$ 
  are Gorenstein rings. Similar results
  hold if one replaces $\Gamma_1(N)$ and $J_1(N)$ by $\Gamma_0(N)$
  and~$J_0(N)$.
\end{prop}

\pf
For the proof we have to pass to $\TT_{\ZZ_p} = \TT_\ZZ \otimes_\ZZ \ZZ_p$.
We also denote by~$\m$ the maximal ideal in~$\TT_{\ZZ_p}$ that corresponds to~$\m$.
Let $V$ be the $\m$-part of the $p$-Tate module of $J_1(N)_\QQ$.
Multiplicity one implies that $V/\m V$ is a $2$-dimensional
$\TT/\fm = \TT_{\ZZ_p}/\m = \TT_{\FF_p,\widetilde{\m}}/\widetilde{\m}$-vector 
space, since
$$ V/\m V \cong (V/pV)/\widetilde{\m}
          \cong (J_1(N)_\QQ(\Qbar)[p]) / \widetilde{\m} 
          \cong (J_1(N)_\QQ(\Qbar)[p])^\vee / \widetilde{\m}
          \cong (J_1(N)_\QQ(\Qbar)[\m])^\vee,$$
where the self-duality comes from the modified Weil pairing which
respects the Hecke action (see e.g.\ \cite{gross}, p.~485).  Nakayama's
Lemma hence implies that $V$ is a $\TT_{\ZZ_p,\m}$-module of
rank~$2$.  As one knows that $V \otimes_{\ZZ_p} \QQ_p$ is a $\TT_\m
\otimes \QQ_p$-module of rank~$2$, it follows that $V$ is a free
$\TT_{\ZZ_p,\m}$-module of rank~$2$, whence
$J_1(N)_\QQ(\Qbar)[p]_{\widetilde{\m}}$ is a free
$\TT_{\FF_p,\widetilde{\m}}$-module of rank~$2$.  Taking the
$\widetilde{\m}$-kernel gives 
$J_1(N)_\QQ(\Qbar)[\m]=(\TT_{\FF_p,\widetilde{\m}}[\widetilde{\m}])^2$,
whence the Gorenstein defect is zero.
In the $\Gamma_0$-situation, the same proof holds.
\qed
\medskip

In the so-called ordinary case, we have the following precise relationship
between the multiplicity and the Gorenstein defect, which was
suggested to us by Kevin Buzzard. The proof can be found in~\cite{nongor}.

\begin{prop}\label{gormult}
  Suppose $p \nmid N$ and let $M=N$ or $M=Np$.  Let $\fm$ be a maximal
  ideal of $\heckeoneNk M2$ of residue characteristic~$p$ and assume
  that $\fm$ is ordinary,
  i.e.\ that the $p$-th Hecke operator~$T_p$ is not in~$\fm$. Assume also
  that $\rho_\m$ is irreducible. Denote
  by $\widetilde{\fm}$ the image of~$\m$ in 
  $\TT_{\FF_p} := \TT_{\ZZ\to\FF_p}(S_2(\Gamma_1(M)))$. 
  Then the following statements hold:
\begin{enumerate}[(a)]
\item 
  There is the exact sequence
  $$ 0 \to \TT_{\FF_p,\widetilde{\fm}} 
       \to J_1(M)(\Qbar)[p]_{\widetilde{\fm}} 
       \to \TT_{\FF_p,\widetilde{\fm}}^\vee 
       \to 0$$ 
  of $\TT_{\FF_p,\widetilde{\fm}}$-modules, where the dual is the
  $\FF_p$-linear dual.

\item 
  If $d$ is the Gorenstein defect of $\TT_{\FF_p,\widetilde{\fm}}$
  and $r$ is the multiplicity of $\rho_\fm$, then the relation
  $$ d = 2r-2$$
  holds.
\end{enumerate}
\end{prop}

We now establish a relation between mod~$p$ Hecke algebras of weights
$3 \le k \le p+1$ for levels~$N$ not divisible by~$p$ with Hecke
algebras of weight~$2$ and level~$Np$. It is needed in order to
compare the Hecke algebras in higher weight to those acting on the
$p$-torsion of Jacobians and thus to make a link to the multiplicity
of the attached Galois representations.

\begin{prop}\label{totwo}
  Let $N \ge 5$, $p \nmid N$ and $3 \le k \le p+1$. Let $\m$ be a maximal
  ideal of the mod~$p$ Hecke algebra $\TT_{\ZZ\to\FF_p}(S_k(\Gamma_1(N))$.
  Then there exists a maximal ideal $\m_2$
  of $\TT_{\ZZ\to\FF_p}(S_2(\Gamma_1(Np))$ and a natural surjection
  $$\TT_{\ZZ\to\FF_p}(S_2(\Gamma_1(Np))_{\m_2} \twoheadrightarrow
    \TT_{\ZZ\to\FF_p}(S_k(\Gamma_1(N))_\m.$$
  If $\m$ is ordinary ($T_p \not\in \m$), this surjection is an
  isomorphism.
\end{prop}

\pf 
From Sections~5 and~6 of \cite{Faithful}, whose notation we adopt
for this proof, one obtains without difficulty the 
commutative diagram of Hecke algebras:
$$ \xymatrix@=1.0cm{
\TT_{\ZZ\to\FF_p}(S_2(\Gamma_1(Np))_{\m_2} \ar@{->>}[r] \ar@{=}[d] & 
\TT_{\FF_p}(J_1(Np)_\QQ(\Qbar)[p])_{\m_2} \ar@{->>}[r] \ar@{->>}[d] & 
\TT_{\FF_p} (H^1(\Gamma_1(N),V_{k-2}(\FF_p)))_\m \\
\TT_{\ZZ_p \to \FF_p} (L)_{\m_2} \ar@{->>}[r] &
\TT_{\FF_p} (\overline{L})_{\m_2} \ar@{->>}[r] &
\TT_{\ZZ\to\FF_p}(S_k(\Gamma_1(N))_\m. \ar@{->>}[u]}$$%
The claimed surjection can be read off. In the ordinary situation,
Proposition~\ref{gormult} shows that the upper left horizontal arrow
is in fact an isomorphism. That also the upper right horizontal arrow
is an isomorphism is explained in~\cite{Faithful}. The result follows
immediately.  
\qed
\medskip

In the next proposition we compare Hecke algebras for spaces of modular
forms on $\Gamma_1(N)$ to those of the same level and weight, but with a
Dirichlet character.

\begin{prop}
  Let $N \ge 5$, $k \ge 2$ and let $\chi: (\ZZ/N\ZZ)^\times \to
  \CC^\times$ be a Dirichlet character.  Let $f \in S_k(N,\chi)
  \subset S_k(\Gamma_1(N))$ be a normalised Hecke eigenform. Let
  further $\m_\chibar$ be the maximal ideal in $\TT_\FF^\chibar :=
  \TT_{\ZZ[\chi]\to\FF}(S_k(N,\chi))$ and $\m$ the one in
  $\TT_{\ZZ\to\FF_p}(S_k(\Gamma_1(N)))$ of residue characteristic~$p$
  for $p \nmid N$ belonging to~$f$. If $k=2$, suppose additionally
  that $\rho_\m$ is irreducible. If $p=2$, suppose that $\rho_\m$ is
  not induced from $\QQ(\sqrt{-1})$, and if $p=3$, suppose that
  $\rho_\m$ is not induced from $\QQ(\sqrt{-3})$.

  Then the Gorenstein defects of $\TT_{\ZZ[\chi]\to\FF}(S_k(N,\chi))_{\m_\chibar}$ 
  and $\TT_{\ZZ\to\FF_p}(S_k(\Gamma_1(N)))_\m$ are equal.
\end{prop}

\pf
Write $\Delta := (\ZZ/N\ZZ)^\times$ and let $\Delta_p$ be its
$p$-Sylow subgroup. Let $\chibar: \Delta \to \FF^\times$ be the
reduction of~$\chi$ obtained by composing~$\chi$
with~$\ZZ[\chi]\to\FF$.  As the Gorenstein defect is invariant under
base extension, it is no loss to work with $\TT_\FF :=
\TT_{\ZZ\to\FF}(S_k(\Gamma_1(N)))$. We still write~$\m$ for the
maximal ideal in $\TT_\FF$ belonging to~$f$.  Note that $\langle
\delta \rangle - \chibar(\delta) \in \m$ for all~$\delta \in \Delta$.

We let $\Delta$ act on $\TT_\FF$ via the diamond operators and we let
$\FF^\chibar$ be a copy of~$\FF$ with $\Delta$-action through the
inverse of~$\chibar$.
We have
$$       (\TT_{\FF,\m} \otimes_\FF \FF^\chibar)_\Delta 
       = (\TT_{\FF,\m} \otimes_\FF \FF^\chibar)/(1-\delta | \delta \in \Delta)
   \cong \TT_{\FF,\m_\chibar}^\chibar,$$
which one obtains by considering the duals, identifying Katz cusp forms
with mod~$p$ ones on~$\Gamma_1(N)$ and applying Carayol's
Lemma (\cite{edixhoven-boston}, Prop.~1.10). For the case $k=2$, we should
point the reader to the correction at the end of the introduction 
to~\cite{edixhoven}. However, the statement still holds after
localisation at maximal ideals corresponding to irreducible representations.
Moreover, the equality
$   (\TT_{\FF,\m} \otimes_\FF \FF^\chibar)^\Delta 
   = \TT_{\FF,\m}[\langle \delta \rangle - \chibar(\delta) | \delta \in \Delta ]$
holds by definition.

Now Lemma~7.3 of~\cite{Faithful} tells us that the localisation at~$\m$
of the $\FF$-vector space of Katz cusp forms of weight~$k$ on $\Gamma_1(N)$ over~$\FF$
is a free $\FF[\Delta_p]$-module. Note that the standing hypothesis $k \le p+1$
of Section~7 of~\cite{Faithful} is not used in the proof of that lemma
and see also \cite{Faithful},~Remark~7.5.
From an elementary calculation one now obtains that 
$N_\Delta = \sum_{\delta \in \Delta} \delta$ induces an isomorphism
$$ (\TT_{\FF_\m} \otimes_\FF \FF^\chibar)_\Delta
   \xrightarrow{N_\Delta}
   (\TT_{\FF_\m} \otimes_\FF \FF^\chibar)^\Delta.$$
We now take the $\m_\chibar$-kernel on both sides and obtain
$$ \TT_{\FF,\m_\chibar}^\chibar [\m_\chibar]
   \cong
   (\TT_{\FF,\m} \otimes_\FF \FF^\chibar)_\Delta [\m_\chibar]
   \cong
   (\TT_{\FF,\m} \otimes_\FF \FF^\chibar)_\Delta [\m]
   \cong
   (\TT_{\FF,\m} \otimes_\FF \FF^\chibar)^\Delta [\m]
   =
   \TT_{\FF,\m}[\m].$$
This proves that the two Gorenstein defects are indeed equal.
\qed
\medskip

The Gorenstein defect that we calculate on the computer is the
number~$d$ of the following corollary, which relates it to the
multiplicity of a Galois representation.

\begin{cor}
  Let $p$ be a prime, $N \ge 5$ an integer such that $p \nmid N$, $k$
  an integer satisfying $2 \le k \le p$ and $\chi: (\ZZ/N\ZZ)^\times
  \to \CC^\times$ a character. Let $f \in S_k(N,\chi)$ be a
  normalised Hecke eigenform. Let further $\m$ denote the maximal
  ideal in $\TT_{\ZZ[\chi]\to\FF}(S_k(N,\chi))$ belonging to~$f$.
  Suppose that $\m$ is ordinary and that $\rho_\m$ is irreducible and
  not induced from $\QQ(\sqrt{-1})$ (if $p=2$) and not induced from
  $\QQ(\sqrt{-3})$ (if $p=3$). We define~$d$ to be the Gorenstein
  defect of $\TT_{\ZZ[\chi]\to\FF}(S_k(N,\chi))_\m$ and~$r$ to be the
  multiplicity of~$\rho_\m$.
 
  Then the equality $d = 2r-2$ holds.
  \qed
\end{cor}

We include the following proposition because it establishes equality
of the two different notions of multiplicities of Galois
representations in the case of the trivial character.

\begin{prop}\label{nulleins}
  Let $N \ge 1$ and $p \nmid N$ and $f \in S_2(\Gamma_0(N))
  \subseteq S_2(\Gamma_1(N))$ be a normalised Hecke eigenform
  belonging to maximal ideals~$\m_0 \subseteq \TT_{\ZZ\to\FF_p}(S_2(\Gamma_0(N)))$ 
  and $\m_1 \subseteq \TT_{\ZZ\to\FF_p}(S_2(\Gamma_1(N)))$ of residue characteristic~$p$.
  Suppose that $\rho_{\m_0} \cong \rho_{\m_1}$ is irreducible.

  Then the multiplicity of $\rho_{\m_1}$ on
  $J_1(N)_\QQ(\Qbar)[p]$ is equal to the multiplicity of~$\rho_{\m_0}$ on
  $J_0(N)_\QQ(\Qbar)[p]$. Thus, if $p>2$, this multiplicity is equal
  to one by Theorem~\ref{gorenstein-theorem}.
\end{prop}

\pf
Let $\Delta := (\ZZ/N\ZZ)^\times$. We first remark that one has the isomorphism
$$ J_0(N)_\QQ(\Qbar)[p]_{\m_0} \cong \big((J_1(N)_\QQ(\Qbar)[p])^\Delta\big)_{\m_0},$$
which one can for example obtain by comparing with the parabolic cohomology
with $\FF_p$-coefficients of the modular curves $Y_0(N)$ and $Y_1(N)$.
Taking the $\m_0$-kernel yields
$$ J_0(N)_\QQ(\Qbar)[\m_0] \cong J_1(N)_\QQ(\Qbar)[\m_1],$$
since $\m_1$ contains $\langle \delta \rangle - 1 $ for all~$\delta \in \Delta$.
\qed

\section{Modular Symbols and Hecke Algebras}
\label{modular-symbols}

The aim of this section is to present the algorithm that we use for
the computations of local factors of Hecke algebras of mod~$p$ modular
forms.  It is based on mod~$p$ modular symbols which have been
implemented in {\sc Magma} \cite{magma} by William Stein.

The bulk of this section deals with proving the main advance, namely a
stop criterion (Corollary~\ref{stopcor}), which in practice greatly
speeds up the computations in comparison with ``standard''
implementations, as it allows us to work with many fewer Hecke
operators than indicated by the theoretical Sturm bound
(Proposition~\ref{propsturm}). We shall list results proving
that the stop criterion is attained in many cases. However, the stop
criterion does not depend on them, in the sense that it being attained
is equivalent to a proof that the algebra it outputs is equal to a
direct factor of a Hecke algebra of mod~$p$ modular forms.

Whereas for Section~\ref{section-relation} the notion of Katz
modular forms seems the right one, the present section works
entirely with mod~$p$ modular forms, the definition of which
is also recalled. This is very natural, since all results in
this section are based on a comparison with the characteristic
zero theory.

\subsection{Mod~$p$ modular forms and modular symbols}

\subsubsection*{Mod~$p$ modular forms}

Let us for the time being fix integers $N \ge 1$ and $k \ge 2$, as
well as a character $\chi: (\ZZ/N\ZZ)^\times \to \CC^\times$ such that
$\chi(-1) = (-1)^k$. Let $M_k(N,\chi)$ be the space of
holomorphic modular forms for $\Gamma_1(N)$, Dirichlet character
$\chi$, and weight~$k$. It decomposes as a direct sum (orthogonal
direct sum with respect to the Petersson inner product) of its
cuspidal subspace $S_k(N,\chi)$ and its Eisenstein subspace
$\Eis_k(N,\chi)$. As before, we let $\cO = \ZZ[\chi]$.
Moreover, we let $\fP$ be a maximal ideal of~$\cO$ above~$p$ 
with residue field~$\FF$ and $\widehat{\cO}$ be the completion 
of~$\cO$ at~$\fP$. Furthermore, let
$K = \QQ_p(\chi)$ be the field of fractions of~$\widehat{\cO}$
and $\bar{\chi}$ be $\chi$ followed by the natural projection
$\cO \twoheadrightarrow \FF$.

Denote by $M_k(N,\chi\,;\,\cO)$ the sub-$\cO$-module generated by
those modular forms whose (standard) $q$-expansion has coefficients
in~$\cO$. It follows from the $q$-expansion principle that
$$M_k(N,\chi\,;\,\cO) \cong 
  \Hom_{\cO} \big(\TT_{\cO}(M_k(N,\chi)),\cO\big)$$ 
and that hence
$M_k(N,\chi\,;\,\cO) \otimes_{\cO} \CC \cong M_k(N,\chi)$.
We put 
$$ M_k(N,\bar{\chi}\,;\,\FF) := 
  M_k(N,\chi\,;\,\cO) \otimes_{\cO} \FF \cong
  \Hom_{\FF} \big(\TT_{\cO}(M_k(N,\chi)),\FF\big) $$
and call the elements of this space {\em mod~$p$ modular forms}.
The Hecke algebra $\TT_{\cO}(M_k(N,\chi))$ acts naturally and it
follows that 
$\TT_{\cO\to\FF}(M_k(N,\chi)) \cong \TT_{\FF}(M_k(N,\chi\,;\,\FF))$.
Similar statements hold for the cuspidal and the Eisenstein subspaces and
we use similar notations.

We call a maximal ideal $\fm$ of 
$\TT_{\cO\to\FF} (M_k(N,\chi\,;\,\cO)$
(respectively, the corresponding maximal ideal of 
$\TT_{\cO\to\widehat{\cO}} (M_k(N,\chi\,;\,\cO))$)
{\em non-Eisenstein} if and only if
$$ S_k(N,\bar{\chi}\,;\,\FF)_\fm \cong M_k(N,\bar{\chi}\,;\,\FF)_\fm.$$
Otherwise, we call $\fm$ {\em Eisenstein}.
\medskip

We now include a short discussion of minimal and maximal primes,
in view of Proposition~\ref{commalg}.
Write $\TT_{\widehat{\cO}}$ for 
$\TT_{\cO \to \widehat{\cO}} (S_k(N,\chi))$.
Let $\fm$ be a maximal ideal of $\TT_{\widehat{\cO}}$. It corresponds
to a $\Gal(\Fpbar/\FF)$-conjugacy class of
normalised eigenforms in $S_k(N,\bar{\chi}\,;\,\FF)$.
That means for each $n\in \NN$ that the minimal polynomial of
$T_n$ acting on $S_k(N,\bar{\chi}\,;\,\FF)_\fm$
is equal to a power of the minimal polynomial of the coefficient~$a_n$
of each member of the conjugacy class.
Moreover, $\fp$ corresponds precisely as above to a 
$\Gal(\Qpbar/K)$-conjugacy
class of normalised eigenforms in $S_k(N,\chi\,;\,\cO) \otimes_\cO K$.  

Suppose that $\fm$ contains minimal primes $\fp_i$ for $i=1,\dots,r$.
Then the normalised eigenforms corresponding to the~$\fp_i$ are
congruent to one another modulo a prime above~$p$. Conversely, any
congruence arises in this way. Thus, a maximal ideal $\fm$ of
$\TT_{\widehat{\cO}}$ is Eisenstein if and only if it contains a minimal
prime corresponding to a conjugacy class of Eisenstein series.
As it is the reduction of a reducible representation,
the mod~$p$ Galois representation corresponding to a non-Eisenstein prime
is reducible. It should be possible to show the converse.

\subsubsection*{Modular symbols}

We now recall the modular symbols formalism and prove two
useful results on base change and torsion. The main
references for the definitions are \cite{SteinBook}
and~\cite{HeckeMS}.

Let $R$ be a ring, $\Gamma \le \SL_2(\ZZ)$ a subgroup and $V$ a left
$R[\Gamma]$-module. Recall that $\PP^1(\QQ) = \QQ \cup \{\infty\}$ is
the set of cusps of $\SL_2(\ZZ)$, which carries a natural
$\SL_2(\ZZ)$-action via fractional linear transformations. We define
the $R$-modules
$$ \cM_R := R[\{\alpha,\beta\}| \alpha,\beta \in \PP^1(\QQ)]/
\langle \{\alpha,\alpha\}, \{\alpha,\beta\} + \{\beta,\gamma\} +
\{\gamma,\alpha\} | \alpha,\beta,\gamma \in \PP^1(\QQ)\rangle$$ 
and $\cB_R := R[\PP^1(\QQ)]$. 
They are connected via the {\em boundary map} 
  $\delta: \cM_R \to \cB_R$ 
which is given by $\{\alpha,\beta\} \mapsto \beta - \alpha.$
Both are equipped with the natural left $\Gamma$-actions. Also let
$\cM_R(V) := \cM_R \otimes_R V$ and $\cB_R(V) := \cB_R \otimes_R V$
with the left diagonal $\Gamma$-action.
We call the $\Gamma$-coinvariants
$$ \cM_R (\Gamma,V) :=  \cM_R(V)_\Gamma = 
\cM_R(V)/ \langle (x - g x) | g \in \Gamma, x \in \cM_R(V) \rangle$$ 
{\em the space of $(\Gamma,V)$-modular symbols.}
Furthermore, {\em the space of $(\Gamma,V)$-boundary symbols}
is defined as the $\Gamma$-coinvariants
$$ \cB_R(\Gamma,V) :=  \cB_R(V)_\Gamma = 
\cB_R(V)/ \langle (x - g x) | g \in \Gamma, x \in \cB_R(V) \rangle.$$ 
The boundary map $\delta$ induces the {\em boundary map} 
$\cM_R(\Gamma,V) \to \cB_R(\Gamma,V)$.
Its kernel is denoted by $\cCM_R(\Gamma,V)$ and is called 
{\em the space of cuspidal $(\Gamma,V)$-modular symbols.}

Let now $N \ge 1$ and $k \ge 2$ be integers and 
  $\chi: (\ZZ/N\ZZ)^\times \to R^\times$ 
be a character, i.e.\ a group homomorphism, such that 
  $\chi(-1) = (-1)^k$ 
in~$R$. Write $V_{k-2}(R)$ for the homogeneous polynomials of 
degree $k-2$ over~$R$ in two variables, equipped with the natural 
$\Gamma_0(N)$-action.
Denote by $V_{k-2}^\chi (R)$ the tensor product 
  $V_{k-2}(R) \otimes_R R^\chi$
for the diagonal $\Gamma_0(N)$-action which on $R^\chi$ comes from
the isomorphism 
  $\Gamma_0(N) / \Gamma_1(N) \cong (\ZZ/N\ZZ)^\times$ 
given by sending $\mat abcd$ to~$d$ followed by~$\chi^{-1}$.

We use the notation $\cM_k(N,\chi\,;\,R)$ for
$\cM(\Gamma_0(N),V_{k-2}^\chi(R))$, as well as similarly for the
boundary and the cuspidal spaces.  The natural action of the matrix
$\eta = \mat {-1}001$ gives an involution on all of these spaces. We
will denote by the superscript ${}^+$ the subspace invariant under
this involution, and by the superscript ${}^-$ the anti-invariant one.
On all modules discussed so far one has Hecke operators $T_n$ for all
$n \in \NN$ and diamond operators. For a definition
see~\cite{SteinBook}.

\begin{lem}\label{basechange}
  Let $R$, $\Gamma$ and $V$ as above and let $R \to S$ be a ring
  homomorphism. Then
  $$\cM(\Gamma,V) \otimes_R S \cong \cM(\Gamma,V \otimes_R S).$$
\end{lem}

\pf 
This follows immediately from the fact that tensoring and taking
coinvariants are both right exact.  
\qed

\begin{prop}\label{torsion}
Let $R$ be a local integral domain of characteristic zero with
principal maximal ideal $\fm = (\pi)$ and residue field $\FF$ of
characteristic~$p$.  Also let $N \ge 1$, $k \ge 2$ be integers and
$\chi: (\ZZ/N\ZZ)^\times \to R^\times$ a character such that
$\chi(-1) = (-1)^k$.  Suppose 
   (i) that $p \ge 5$ or 
  (ii) that $p =  2$ and $N$ is divisible by a prime which is 
       $3$ modulo~$4$ or by~$4$ or 
 (iii) that $p = 3$ and $N$ is divisible by a prime which is $2$
       modulo~$3$ or by~$9$.  
Then the following statements hold:
\begin{enumerate}[(a)]
\item If $k \ge 3$, then 
$\cM_k(N,\chi\,;\,R)[\pi] = \big(V_{k-2}^\chi(\FF)\big)^{\Gamma_0(N)}$.
\item If $k=2$ or if $3 \le k \le p+2$ and $p \nmid N$, then 
$\cM_k(N,\chi\,;\,R)[\pi] = 0$.
\end{enumerate}
\end{prop}

\pf
The conditions assure that the group $\Gamma_0(N)$ does not have
any stabiliser of order $2p$ for its action on the upper half plane.
Hence, by \cite{HeckeMS}, Theorem~6.1, the modular symbols space
$\cM_k(N,\chi\,;\,R)$ is isomorphic to $H^1(\Gamma_0(N),V_{k-2}^\chi(R))$.
The arguments are now precisely those of the beginning
of the proof of \cite{Faithful}, Proposition~2.6.
\qed

\subsubsection*{Hecke algebras of modular symbols and the Eichler-Shimura isomorphism}

From Lemma~\ref{basechange} one deduces a natural surjection
\begin{equation}\label{eqms}
  \TT_{\cO \to \FF}       (\cM_k(N,\chi\,;\,\cO)) 
    \twoheadrightarrow
  \TT_{\FF}(\cM_k(N,\bar{\chi}\,;\,\FF)).
\end{equation}
In the same way, one also obtains 
\begin{equation}\label{eqmseins}
   \TT_{\cO}    (\cM_k(N,\chi\,;\,\cO)) \twoheadrightarrow
   \TT_{\cO}    (\cM_k(N,\chi\,;\,\cO)/\text{torsion})
   \cong \TT_{\cO}    (\cM_k(N,\chi\,;\,\CC)),
\end{equation}
where one uses for the isomorphism that the Hecke operators are already
defined over~$\cO$. A similar statement holds for the cuspidal subspace.

We call a maximal prime $\fm$ of 
$\TT_{\cO\to\widehat{\cO}} (\cM_k(N,\chi\,;\,\cO))$
(respectively the corresponding prime of 
$\TT_{\cO\to\FF} (\cM_k(N,\chi\,;\,\cO))$)
{\em non-torsion} if
$$  \cM_k(N,\chi\,;\,\widehat{\cO})_\fm \cong 
   (\cM_k(N,\chi\,;\,\widehat{\cO})/\text{torsion})_\fm.$$
This is equivalent to the height of $\fm$ being~$1$.
Proposition~\ref{torsion} tells us some cases in which all
primes are non-torsion.

\begin{thm}[Eichler-Shimura]\label{thmes}
There are isomorphisms respecting the Hecke operators
\begin{enumerate}[(a)]
\item $M_k(N,\chi) \oplus S_k(N,\chi)^\vee
  \cong \cM_k(N,\chi\,;\,\CC),$
\item $S_k(N,\chi) \oplus S_k(N,\chi)^\vee 
  \cong \cCM_k(N,\chi\,;\,\CC),$
\item $S_k(N,\chi) \cong \cCM_k(N,\chi\,;\,\CC)^+.$
\end{enumerate}
\end{thm}

\pf
Parts~(a) and (b) are \cite{DiamondIm}, Theorem 12.2.2, together with
the comparison of \cite{HeckeMS}, Theorem~6.1. 
We use that the space of anti-holomorphic cusp forms is dual
to the space of holomorphic cusp forms.
Part~(c) is a direct consequence of~(b).
\qed

\begin{cor}\label{cores}
There are isomorphisms
$$\TT_{\cO}(S_k(N,\chi)) \cong \TT_{\cO} (\cCM_k(N,\chi\,;\,\CC))
 \cong \TT_{\cO} (\cCM_k(N,\chi\,;\,\CC)^+),$$
given by sending $T_n$ to $T_n$ for all positive~$n$.
\qed
\end{cor}

\subsection{The stop criterion}

Although it is impossible to determine a priori the dimension of the
local factor of the Hecke algebra associated with a given modular form
mod~$p$, Corollary~\ref{stopcor} implies that the computation of Hecke
operators can be stopped when the algebra generated has reached a
certain dimension that is computed along the way.  This criterion has
turned out to be extremely useful and has made possible some of our
computations which would not have been feasible using the Hecke bound
naively. See Section~\ref{computational-results} for a short
discussion of this issue.

\subsubsection*{Some commutative algebra}

We collect some useful statements from commutative algebra,
which will be applied to Hecke algebras in the sequel.

\begin{prop}\label{commalg}
Let $R$ be an integral domain of characteristic zero which is a finitely
generated $\ZZ$-module. Write $\widehat{R}$ for the completion of $R$ at a 
maximal ideal of~$R$ and denote by $\FF$ the residue field
and by $K$ the fraction field of~$\widehat{R}$.
Let furthermore $A$ be a commutative $R$-algebra which is finitely
generated as an $R$-module.
For any ring homomorphism $R\to S$ write $A_S$ for $A \otimes_R S$.
Then the following statements hold.

\begin{enumerate}[(a)]
\item The Krull dimension of $A_{\widehat{R}}$ is less than or equal
  to~$1$. The maximal ideals of $A_{\widehat{R}}$ correspond
  bijectively under taking pre-images to the maximal ideals of
  $A_\FF$. Primes $\fp$ of height $0$ which are contained in a prime
  of height $1$ of $A_{\widehat{R}}$ are in bijection with primes of
  $A_K$ under extension (i.e.\ $\fp A_K$), for which the notation
  $\fp^e$ will be used.

  Under these correspondences, one has
  $A_{\FF,\fm} \cong A_{\widehat{R},\fm} \otimes_{\widehat{R}} \FF$
  and
  $A_{K,\fp^e} \cong A_{\widehat{R},\fp}$.

\item The algebra $A_{\widehat{R}}$ decomposes as
  $$ A_{\widehat{R}} \cong \prod_\fm A_{\widehat{R},\fm},$$
  where the product runs over the maximal ideals $\fm$ of
  $A_{\widehat{R}}$.

\item The algebra $A_\FF$ decomposes as
  $$ A_\FF \cong \prod_\fm A_{\FF,\fm},$$
  where the product runs over the maximal ideals $\fm$ of $A_\FF$.

\item The algebra $A_K$ decomposes as
  $$ A_K \cong \prod_\fp A_{K,\fp^e} \cong  \prod_\fp A_{\widehat{R},\fp},$$
  where the products run over the minimal prime ideals $\fp$ of
  $A_{\widehat{R}}$ which are contained in a prime ideal of
  height~$1$.
\end{enumerate}
\end{prop}

\pf 
As $A_{\widehat{R}}$ is a finitely generated $\widehat{R}$-module,
$A_{\widehat{R}}/\fp$ with a prime $\fp$ is an integral domain which
is a finitely generated $\widehat{R}$-module. Hence, it is either a
finite field or a finite extension of $\widehat{R}$. This proves that
the height of $\fp$ is less than or equal to~$1$. The correspondences
and the isomorphisms of Part~(a) are easily verified. The
decompositions in Parts~(b) and~(c) are \cite{Eisenbud},
Corollary~7.6.  Part~(d) follows by tensoring~(b) over $\widehat{R}$
with~$K$.
\qed
\medskip

Similar decompositions for $A$-modules are derived by applying the
idempotents of the decompositions of Part~(b).

\begin{prop}\label{dimension}
Assume the set-up of Proposition~\ref{commalg} and
let $M,N$ be $A$-modules which as $R$-modules are free of finite rank.
Suppose that
\begin{enumerate}[(a)]
\item $M \otimes_R \CC \cong N \otimes_R \CC$
as $A \otimes_R \CC$-modules, or
\item $M \otimes_R \bar{K} \cong N \otimes_R \bar{K}$
as $A \otimes_R \bar{K}$-modules.
\end{enumerate}
Then for all prime ideals $\fm$ of $A_\FF$ corresponding
to height $1$ primes of $A_{\widehat{R}}$ the equality
$$ \dim_\FF (M \otimes_R \FF)_\fm = \dim_\FF (N \otimes_R \FF)_\fm$$
holds.
\end{prop}

\pf
As for $A$, we also write $M_K$ for $M \otimes_R K$ and similarly for
$N$ and $\widehat{R}$, $\FF$, etc.  By choosing an isomorphism $\CC \cong
\bar{K}$, it suffices to prove Part~(b). 
Using Proposition~\ref{commalg}, Part~(d), the isomorphism
$M \otimes_R \bar{K} \cong N \otimes_R \bar{K}$ can be rewritten as
$$ \bigoplus_{\fp} (M_{K, \fp^e} \otimes_K \bar{K}) \cong
\bigoplus_{\fp} (N_{K, \fp^e} \otimes_K \bar{K}),$$
where the sums run over the minimal primes $\fp$ of $A_{\widehat{R}}$
which are properly contained in a maximal prime.  Hence, an
isomorphism
$M_{K, \fp^e} \otimes_K \bar{K} \cong N_{K, \fp^e} \otimes_K \bar{K}$
exists for each~$\fp$.  Since for each maximal ideal $\fm$ of
$A_{\widehat{R}}$ of height $1$ we have by Proposition~\ref{commalg}
$$ M_{\widehat{R}, \fm} \otimes_{\widehat{R}} K \cong 
   \bigoplus_{\fp \subseteq \fm \text{ min.}} M_{K, \fp^e}$$ 
and similarly for $N$, we get
\begin{align*}  
    \dim_\FF M_{\FF, \fm} 
  =& \rk_{\widehat{R}} M_{\widehat{R}, \fm} 
  = \sum_{\fp \subseteq \fm \text{ min.}} \dim_K M_{K, \fp^e} \\
  =& \sum_{\fp \subseteq \fm \text{ min.}} \dim_K N_{K, \fp^e}
  =  \rk_{\widehat{R}} N_{\widehat{R}, \fm} 
  = \dim_\FF N_{\FF, \fm}.
\end{align*}
This proves the proposition.
\qed

\subsubsection*{The stop criterion}

\begin{prop}\label{propdim}
Let $\fm$ be a maximal ideal of
$\TT_{\cO\to\FF}(\cM_k(N,\chi\,;\,\cO))$
which is non-torsion and non-Eisenstein. 
Then the following statements hold:
\begin{enumerate}[(a)]
\item $\cCM_k(N,\bar{\chi}\,;\,\FF)_\fm \cong \cM_k(N,\bar{\chi}\,;\,\FF)_\fm$.
\item $2 \cdot \dim_{\FF} S_k(N,\bar{\chi}\,;\,\FF)_\fm =
         \dim_{\FF} \cCM_k(N,\bar{\chi}\,;\,\FF)_\fm$.
\item If $p \neq 2$, then
        $\dim_{\FF} S_k(N,\bar{\chi}\,;\,\FF)_\fm =
         \dim_{\FF} \cCM_k(N,\bar{\chi}\,;\,\FF)^+_\fm$.
\end{enumerate}
\end{prop}

\pf
Part~(c) follows directly from Part~(b) by decomposing 
$\cCM_k(N,\bar{\chi}\,;\,\FF)$
into a direct sum of its plus- and its minus-part.
Statements~(a) and~(b) will be concluded from Proposition~\ref{dimension}.
More precisely, it allows us to derive from Theorem~\ref{thmes} that
\begin{align*}
  &\dim_{\FF}\big((\cM_k(N,\chi\,;\,\cO)/\text{torsion}) 
                                \otimes_{\cO} \FF\big)_\fm\\
= &\dim_{\FF}\big(\Eis_k (N,\bar{\chi}\,;\,\FF) 
                         \oplus S_k (N,\bar{\chi}\,;\,\FF) 
                         \oplus S_k (N,\bar{\chi}\,;\,\FF)^\vee  \big)_\fm
\end{align*}
and 
$$ \dim_{\FF}\big((\cCM_k(N,\chi\,;\,\cO)/\text{torsion}) 
                              \otimes_{\cO} \FF\big)_\fm
=  2 \cdot \dim_{\FF}\big(S_k(N,\bar{\chi}\,;\,\FF)\big)_\fm.$$
The latter proves Part~(b), since $\fm$ is non-torsion.
As by the definition of a non-Eisenstein prime 
$\Eis_k (N,\bar{\chi}\,;\,\FF)_\fm = 0$ and
again since $\fm$ is non-torsion, it follows that
$$\dim_{\FF}\cCM_k(N,\bar{\chi}\,;\,\FF)_\fm 
= \dim_{\FF}\cM_k(N,\bar{\chi}\,;\,\FF)_\fm,$$
which implies Part~(a).
\qed
\medskip

We will henceforth often regard non-Eisenstein non-torsion primes 
as in the proposition as maximal primes of 
$\TT_{\FF}(S_k(N,\bar{\chi}\,;\,\FF)) = \TT_{\cO\to\FF}(S_k(N,\chi))$.

\begin{cor}[Stop Criterion]\label{stopcor}
Let $\fm$ be a maximal ideal of $\TT_{\FF}(S_k(N,\bar{\chi}\,;\,\FF))$
which is non-Eisen\-stein and non-torsion.
\begin{enumerate}[(a)]
\item One has $\dim_{\FF} \cM_k(N,\bar{\chi}\,;\,\FF)_\fm = 
2 \cdot \dim_{\FF} \TT_{\FF} \big(\cM_k(N,\bar{\chi}\,;\,\FF)\big)_\fm$
if and only if 
$$\TT_{\FF} \big(S_k(N,\bar{\chi}\,;\,\FF)\big)_\fm \cong
\TT_{\FF} \big(\cCM_k(N,\bar{\chi}\,;\,\FF)\big)_\fm.$$
\item One has $\dim_{\FF} \cCM_k(N,\bar{\chi}\,;\,\FF)_\fm = 
2 \cdot \dim_{\FF} \TT_{\FF} \big(\cCM_k(N,\bar{\chi}\,;\,\FF)\big)_\fm$
if and only if 
$$\TT_{\FF} \big(S_k(N,\bar{\chi}\,;\,\FF)\big)_\fm \cong
\TT_{\FF} \big(\cCM_k(N,\bar{\chi}\,;\,\FF)\big)_\fm.$$
\item Assume $p \neq 2$. One has $\dim_{\FF} \cCM_k(N,\bar{\chi}\,;\,\FF)^+_\fm = 
\dim_{\FF} \TT_{\FF} \big(\cCM_k(N,\bar{\chi}\,;\,\FF)\big)_\fm$
if and only if 
$$\TT_{\FF} \big(S_k(N,\bar{\chi}\,;\,\FF)\big)_\fm \cong
\TT_{\FF} \big(\cCM_k(N,\bar{\chi}\,;\,\FF)^+\big)_\fm.$$
\end{enumerate}
\end{cor}

\pf
We only prove (a), as (b) and (c) are similar.
From Part~(b) of Proposition~\ref{propdim}
and the fact that the $\FF$-dimension of the algebra
$\TT_{\FF} \big(S_k(N,\bar{\chi}\,;\,\FF)\big)_\fm$
is equal to the one of
$S_k(N,\bar{\chi}\,;\,\FF)$, as they are dual to each other,
it follows that 
$$ 2 \cdot \dim_\FF \TT_{\FF} \big(S_k(N,\bar{\chi}\,;\,\FF)\big)_\fm =
\dim_\FF \big(\cCM_k(N,\bar{\chi}\,;\,\FF)\big)_\fm.$$
The result is now a direct consequence of Equations~\ref{eqms}
and~\ref{eqmseins} and Corollary~\ref{cores}.
\qed
\medskip

Note that the first line of each statement only uses modular symbols
and not modular forms, but it allows us to make statements involving
modular forms.  This is the aforementioned stop criterion; the computation of Hecke
operators can be stopped if this equality is reached.

We now list some results concerning the validity of the equivalent
statements of Corollary~\ref{stopcor}. The reader can also consult
\cite{EPW} for general results in the ordinary and distinguished case.

\begin{prop}
Let $p \ge 5$ be a prime, $k \ge 2$ and $N \ge 5$ with $p \nmid N$ integers,
$\FF$ a finite extension of~$\FF_p$,
$\bar{\chi}: (\ZZ/N\ZZ)^\times \to \FF^\times$ a character
and $\fm$ a maximal ideal of
$\TT_\FF(S_k(N,\bar{\chi}\,;\,\FF))$ which is non-Eisenstein and non-torsion.
Suppose (i) that $2 \le k \le p-1$ or
(ii) that $k \in \{p,p+1\}$ and $\fm$ is ordinary.
Then
$$ \TT_{\FF} \big(S_k(N,\bar{\chi}\,;\,\FF)\big)_\fm  \cong
   \TT_{\FF} \big(\cCM_k(N,\bar{\chi}\,;\,\FF)\big)_\fm \cong
   \TT_{\FF} \big(\cCM_k(N,\bar{\chi}\,;\,\FF)^+\big)_\fm.$$
\end{prop}

\pf
Using the comparison with group cohomology of \cite{HeckeMS}, Theorem~6.1,
the result follows under Assumption~(i) from \cite{edixhoven}, Theorem~5.2,
and is proved under Assumption~(ii) in \cite{Faithful}, Corollary~6.9,
for the case of the group $\Gamma_1(N)$ and no Dirichlet character.
The passage to a character is established by \cite{Faithful}, Theorem~7.4
and the remark following it. One identifies the mod~$p$ modular
forms appearing with corresponding Katz forms using Carayol's Lemma
(\cite{edixhoven-boston}, Prop.~1.10).
\qed
\medskip

We end this section by stating the so-called Sturm bound (also called
the Hecke bound), which gives the best a priori upper bound for how many
Hecke operators are needed to generate all the Hecke algebra. We only
need it in our algorithm in cases in which it is theoretically not
known that the stop criterion will be reached. This will enable the
algorithm to detect if the Hecke algebra on modular symbols is not
isomorphic to the corresponding one on cuspidal modular forms.

\begin{prop}[Sturm bound]\label{propsturm}
The Hecke algebra $\TT_\CC(S_k(N,\chi))$ can be generated
as an algebra by the Hecke operators $T_l$ for all primes~$l$
smaller than or equal to
$\frac{kN}{12}\prod_{q \mid N, q \text{ prime}} (1+\frac{1}{q})$.
\end{prop}

\pf 
This is discussed in detail in Chapter~11 of \cite{SteinBook}.
\qed

\subsection{Algorithm}

In this section we present a sketch of the algorithm that we used for
our computations. The {\sc Magma} code can be downloaded from the
second author's webpage and a manual is included as 
Appendix~\ref{section-manual}.

\noindent
{\bf\underline{Input:}}   Integers $N \ge 1$, $k \ge 2$, a finite field $\FF$,
a character $\chi: (\ZZ/N\ZZ)^\times \to \FF^\times$ and for each prime~$l$
less than or equal to the Sturm bound an irreducible polynomial $f_l \in \FF[X]$.\\
{\bf \underline{Output:}}  An $\FF$-algebra.
\begin{itemize}
\itemsep=0cm plus 0pt minus 0pt
\item $M \leftarrow \cCM_k(N,\chi\,;\,\FF)$, $l \leftarrow 1$,
$L \leftarrow $ empty list.
\item repeat
\begin{itemize}
\itemsep=0cm plus 0pt minus 0pt
\item $l \leftarrow $ next prime after $l$.
\item Compute $T_l$ on $M$ and append it to the list $L$.
\item $M \leftarrow $ the restriction of $M$ to the $f_l$-primary
subspace for $T_l$, i.e.\ to the biggest subspace of $M$ on which the minimal polynomial 
of $T_l$ is a power of $f_l$.
\item $A \leftarrow $ the $\FF$-algebra generated by the restrictions to $M$
of $T_2, T_3, \dots, T_l$.
\end{itemize}
\item until $2 \cdot \dim (A) = \dim (M)$ \emph{[the stop criterion]} or $l > $ Sturm bound.
\item return $A$.
\end{itemize}

The $f_l$ should, of course, be chosen as the minimal polynomials of
the coefficients $a_l(f)$ of the normalised eigenform $f \in
S_k(N,\chi\,;\,\bar{\FF})$ whose local Hecke algebra one wants to
compute.  Suppose the algorithm stops at the prime~$q$.  If $q$ is
bigger than the Sturm bound, the equivalent conditions of
Corollary~\ref{stopcor} do not hold.  In that case the output should
be disregarded.  Otherwise, $A$ is isomorphic to a direct product of
the form $\prod_\fm \TT(S_k(N,\chi\,;\,\FF))_\fm$ where the $\fm$ are
those maximal ideals such that the minimal polynomials of $T_2, T_3,
\dots, T_q$ on $\TT(S_k(N,\chi\,;\,\FF))_\fm$ are equal to powers of
$f_2, f_3, \dots, f_q$.  It can happen that $A$ consists of more than
one factor.  Hence, one should still decompose $A$ into its local
factors.  Alternatively, one can also replace the last line but one in
the algorithm by
\begin{itemize}
\item until $\big((2 \cdot \dim (A) = \dim (M))$ and $A$ is local$\big)$ 
or $l > $ Sturm bound,
\end{itemize}
which ensures that the output is a local algebra.  In practice, one
modifies the algorithm such that not for every prime~$l$ a
polynomial~$f_l$ need be given, but that the algorithm takes each
irreducible factor of the minimal polynomial of~$T_l$ if no $f_l$ is
known. It is also useful to choose the order how $l$ runs through
the primes. For example, one might want to take $l=p$ at an early
stage with $p$ the characteristic of~$\FF$, if one knows that
this operator is needed, which is the case in all computations
concerning Question~\ref{ourquestion}.

\section{Computational results}
\label{computational-results}

In view of Question~\ref{ourquestion}, we produced 384 examples of
odd irreducible continuous Galois representations $\Gal(\Qbar/\QQ) \to
\GL_2(\Fpbar)$ that are completely split at~$p$. The results are
documented in the accompanying tables (Appendix~\ref{section-tables}).
The complete data can be downloaded from the second author's webpage.

The Galois representations were created either by class field theory
or from an irreducible integer polynomial whose Galois group embeds
into $\GL_2(\Fpbar)$. All examples but one are dihedral; the remaining
one is icosahedral. For each of these an eigenform was computed giving
rise to it. The Gorenstein defect of the corresponding local Hecke
algebra factor turned out always to be~$2$, supporting 
Question~\ref{ourquestion}.

The authors preferred to proceed like this, instead of computing all
Hecke algebras mod $p$ in weight $p$ for all ``small'' primes~$p$ and
all ``small'' levels, since non-dihedral examples in which the
assumptions of Question~\ref{ourquestion} are satisfied are very rare.

\subsection{Table entries}

For every computed local Hecke algebra enough data are stored to
recreate it as an abstract algebra and important characteristics
are listed in the tables of Appendix~\ref{section-tables}.
A sample table entry is the following.
\begin{longtable}{||c|c|c|c|c|c|c|c|c|c||}
\hline
Level & Wt & ResD & Dim & EmbDim & NilO & GorDef & \#Ops & \#(p$<$HB) & Gp \\
\hline
5939 & 5 & 3 & 12 & 3 & 5 & 2 & 5 & 366 & $D_{7}$ \\
\hline
\end{longtable}

Each entry corresponds to the Galois conjugacy class of an eigenform
$f$ mod~$p$ with associated local Hecke algebra~$A$. The first and the
second column indicate the level and the weight of~$f$. The latter is
in all examples equal to the characteristic of the base field $k$ (a
finite extention of~$\FF_p$) of the algebra.  Let $\fm_A$ denote the
maximal ideal of~$A$. Then ResD stands for the degree of $K = A/\fm_A$
over~$\FF_p$. Let us consider $A \otimes_k K$. It decomposes into a
direct product of a local $K$-algebra $B$ and its
$\Gal(K/k)$-conjugates. The $K$-dimension of $B$ (which is equal to
the $k$-dimension of~$A$) is recorded in the fourth column.

Let $\fm_B$ be the maximal ideal of~$B$. The {\em embedding dimension}
EmbDim is the $K$-dimension of $\fm_B/\fm_B^2$. By Nakayama's Lemma
this is the minimal number of $B$-generators for~$\fm_B$. 
The {\em nilpotency order} NilO is the maximal integer~$n$ such that
$\fm_B^n$ is not the zero ideal. The column GorDef contains the
Gorenstein defect of~$B$ (which is the same as the Gorenstein defect
of~$A$).

By \#Ops it is indicated how many Hecke operators were used to
generate the algebra~$A$, applying the stop criterion
(Corollary~\ref{stopcor}).  This is contrasted with the number of
primes smaller than the Sturm bound (Proposition~\ref{propsturm}, it
is also called the Hecke bound), denoted by \#(p$<$HB).  One
immediately observes that the stop criterion is very efficient.
Whereas the Sturm bound is roughly linear in the level, in 365 of the
384 calculated examples, less than 10 Hecke operators sufficed, in 252
examples even 5 were enough.

The final column contains the image of the mod~$p$ Galois
representation attached to~$f$ as an abstract group.

\subsection{Dihedral examples}

All Hecke algebras except one in our tables correspond to eigenforms
whose Galois representations are dihedral, since these are by far
the easiest to obtain explicitly as one can use class field theory.
This is explained now.

Let $p$ be a prime and $d$ a square-free integer which is $1$ mod~$4$
and not divisible by~$p$. We denote by $K$ the quadratic
field~$\QQ(\sqrt{d})$.  Further we consider an unramified character
$\chi: \Gal(\Qbar/K) \to \Fpbar^\times$ of order $n \ge 3$. We assume
that its inverse $\chi^{-1}$ is equal to $\chi$ conjugated by
$\sigma$, denoted $\chi^\sigma$, for $\sigma$ (a lift of) the
non-trivial element of $\Gal(K/\QQ)$. The induced representation
$$\rho_\chi := \Ind_{\Gal(\Qbar/K)}^{\Gal(\Qbar/\QQ)}(\chi) : 
  \Gal(\Qbar/\QQ) \to \GL_2(\Fpbar)$$ 
is irreducible and its image is the dihedral group $D_n$ of order~$2n$.
If $l$ is a prime not dividing $2d$, we have $\rho_\chi(\Frob_l) =
\mat 0110$ if $\big(\frac{d}{l}\big)=-1$, and $\rho_\chi(\Frob_l) =
\mat {\chi(\Frob_\Lambda)}00 {\chi^\sigma(\Frob_\Lambda)}$ if
$\big(\frac{d}{l}\big)=1$ and $l \cO_K = \Lambda \sigma(\Lambda)$.
This explicit description makes it obvious that the determinant of
$\rho_\chi$ is the Legendre symbol $l \mapsto \big(\frac{d}{l}\big)$.

Since the kernel of $\chi$ corresponds to a subfield of the Hilbert
class field of~$K$, simple computations in the class group of~$K$
allow one to determine which primes split completely. These give
examples satisfying the assumptions of Question~\ref{ourquestion}
(the Frobenius at~$p$ is the identity)
if $\rho_\chi$ is odd, i.e.\ if $p=2$ or $d < 0$.

We remark that for characters $\chi$ of odd order~$n$ the assumption
$\chi^{-1} = \chi^\sigma$ is not a big restriction, since any character
can be written as $\chi = \chi_1 \chi_2$ with $\chi_1^\sigma = \chi_1^{-1}$
and $\chi_2^\sigma = \chi_2$, hence the latter descends to a character
of $\Gal(\Qbar/\QQ)$ and the representation $\rho_\chi$ is isomorphic
to $\rho_{\chi_1} \otimes \chi_2$.

All dihedral representations are known to come from eigenforms in the
minimal possible weight with level equal to the (outside of $p$) 
conductor of the representation (see \cite{Dihedral}, Theorem~1).
\medskip

In the tables we computed the Hecke algebras of odd dihedral
representations as above in the following ranges. For each prime~$p$
less than~$100$ and each prime $l$ less than or equal to the largest
level occuring in the table for~$p$, we chose $d$ as plus or minus~$l$
such that $d$ is $1$ mod $4$ and we let $H$ run through all
non-trivial cyclic quotients of the class group of $\QQ(\sqrt{d})$ of
order coprime to~$p$. For each $H$ we chose (unramified)
characters~$\chi$ of the absolute Galois group of $\QQ(\sqrt{d})$
corresponding to~$H$, up to Galois conjugacy and up to replacing
$\chi$ by its inverse. Then $\chi$ is not the restriction of a
character of $\Gal(\Qbar/\QQ)$. By genus theory the order of $\chi$ is
odd, as the class number is, so we necessarily have $\chi^{-1} =
\chi^\sigma$.  We computed the local factor of
$\TT_{\FF_p}(S_p(d,\big(\frac{d}{\cdot}\big) \,;\, \FF_p))$ corresponding
to~$\rho_\chi$ if $\rho_\chi$ is odd and $p$ is completely split. For
the prime $p=2$ we also allowed square-free integers~$d$ which are $1$
mod~$4$ and whose absolute value is less than~$5000$.

\subsection{Icosahedral example}

With the help of a list of polynomials provided by Gunter Malle
(\cite{malle}) a Galois representation of $\Gal(\Qbar/\QQ)$ with
values in $\GL_2(\Ftbar)$ which is of prime conductor, completely
split at~$2$ and thus satisfies the assumptions of
Question~\ref{ourquestion} and whose image is isomorphic to the
icosahedral group~$A_5$ could be described explicitly.  The modular
forms in weight $2$ predicted by Serre's conjecture were found and the
corresponding Hecke algebra turned out to have Gorenstein defect equal
to~$2$.

Let $f \in \ZZ[X]$ be an irreducible polynomial of degree~$5$ whose
Galois group, i.e.\ the Galois group of the normal closure $L$ of $K =
\QQ[X]/(f)$, is isomorphic to~$A_5$. We assume that $K$ is unramified
at $2$, $3$ and~$5$. We have the Galois representation
$$ \rho_f : \Gal(\Qbar/\QQ) \twoheadrightarrow 
            \Gal(L/\QQ) \cong A_5 \cong \SL_2(\FF_4).$$ 
We now determine its conductor and its traces.  Let $p$ be a ramified
prime. As the ramification is tame, the image of the inertia group
$\rho_f(I_p)$ at $p$ is cyclic of order $2$, $3$ or~$5$.  In the first
case, the image of a decomposition group $\rho_f(D_p)$ at $p$ is
either equal to $\rho_f(I_p)$ or equal to $\ZZ/2\ZZ \times
\rho_f(I_p)$.  If the order of $\rho_f(I_p)$ is odd and $\rho_f(I_p) =
\rho_f(D_p)$, then any completion of $L$ at the unique prime above~$p$
is totally ramified and cyclic of degree $\#\rho_f(I_p)$, hence
contained in $\QQ_p(\zeta_p)$ for $\zeta_p$ a primitive $p$-th root of
unity. It follows that $p$ is congruent to $1$ mod $\#\rho_f(I_p)$.
If the order of $\rho_f(I_p)$ is odd, but $\rho_f(I_p)$ is not equal
to $\rho_f(D_p)$, then $\rho_f(D_p)$ is a dihedral group and the
completion of $L$ at a prime above $p$ has a unique unramified
quadratic subfield~$S$.  Thus, we have the exact sequence
$$ 0 \to \rho_f(I_p) \to \rho_f(D_p) \to \Gal(S/\QQ_p) \to 0.$$
On the one hand, it is well-known that the conjugation by a lift of the
Frobenius element of $\Gal(S/\QQ_p)$ acts on $\rho_f(I_p)$ by raising
to the $p$-th power. On the other hand, as the action is non-trivial
it also corresponds to inversion on $\rho_f(I_p)$, since the only elements of
order~$2$ in $(\ZZ/3\ZZ)^\times$ and $(\ZZ/5\ZZ)^\times$ are~$-1$.
As a consequence, $p$ is congruent to $-1$ mod $\#\rho_f(I_p)$ in this case.

We hence have the following cases.
\begin{enumerate}[(1)]
\item Suppose $p \cO_K = \fP^5$. Then $p \equiv \pm 1 \mod 5$.
\begin{enumerate}[(a)]
\item If $p \equiv 1 \mod 5$, then $\rho_f|_{I_p} \sim \mat \chi00{\chi^{-1}}$
with $\chi$ a totally ramified character of $\Gal(\Qpbar/\QQ_p)$ of order~$5$.
\item If $p \equiv -1 \mod 5$, then $\rho_f(D_p)$ is the dihedral group with $10$ elements.
\end{enumerate}

\item Suppose $p \cO_K = \fP^3 \fQ \fR$ or $p \cO_K = \fP^3 \fQ$. 
\begin{enumerate}[(a)]
\item If $p \equiv 1 \mod 3$, then $\rho_f|_{I_p} \sim \mat \chi00{\chi^{-1}}$
with $\chi$ a totally ramified character of $\Gal(\Qpbar/\QQ_p)$ of order~$3$.
\item If $p \equiv -1 \mod 3$, then $\rho_f(D_p)$ is the dihedral group with $6$ elements.
\end{enumerate}

\item Suppose that $p$ is ramified, but that we are neither in Case~(1)
nor in Case~(2). Then $\rho_f|_{I_p} \sim \mat 1101$.
\end{enumerate}

By the definition of the conductor at $p$ it is clear that it is $p^2$
in Cases (1) and~(2) and $p$ in Case~(3). However, in Cases (1)(a) and
(2)(a) one can choose a character $\epsilon$ of $\Gal(\Qbar/\QQ)$ of
the same order as~$\chi$ whose restriction to $D_p$ gives the
character~$\chi$. If one twists the representation $\rho_f$ by
$\epsilon$ one finds also in these cases that the conductor at~$p$
is~$p$.

An inspection of the conjugacy classes of the group
$\SL_2(\FF_4)$ shows that the traces of $\rho_f$ twisted
by some character $\epsilon$ of $\Gal(\Qbar/\QQ)$ are as follows.
Let $l$ be an unramified prime.
\begin{itemize}
\item If the order of $\Frob_l$ is $5$, then the trace at $\Frob_l$
is~$\epsilon(\Frob_l) w$ where $w$ is a root of the polynomial $X^2+X+1$
in $\FF_2[X]$.
\item If the order of $\Frob_l$ is $3$, then the trace at $\Frob_l$
is~$\epsilon(\Frob_l)$.
\item If the order of $\Frob_l$ is $1$ or $2$, then the trace at $\Frob_l$
is~$0$.
\end{itemize}

These statements allow the easy identification of the modular form
belonging to an icosahedral representation.

We end this section with some remarks on our icosahedral example.
It was obtained using the polynomial 
$x^5 - x^4 - 79 x^3 + 225 x^2 + 998 x - 3272$.
The corresponding table entry is:
\begin{longtable}{||c|c|c|c|c|c|c|c|c|c||}
\hline
Level & Wt & ResD & Dim & EmbDim & NilO & GorDef & \#Ops & \#(p$<$HB) & Gp \\
\hline
89491 & 2 & 2 & 12 & 4 & 3 & 2 & 4 & 1746 & $A_5$ \\
\hline
\end{longtable}

Hence, in level $89491$ and weight~$2$ there is a single eigenform~$g$
mod~$2$ up to Galois conjugacy whose first couple of $q$-coefficients
agree with the traces of a twist of the given icosahedral Galois
representation. From this one can deduce that the Galois
representation $\rho_g$ of~$g$ has an icosahedral image and is only
ramified at~$89491$. As weight and level lowering are not known in our
case, we cannot prove that $\rho_g$ coincides with a twist of the
given one. It might, however, be possible to exclude the existence of
two distinct icosahedral extensions of the rationals inside~$\CC$ that
ramify only at~$89491$ by consulting tables.  According to Malle, the
icosahedral extension used has smallest discriminant among all totally
real $A_5$-extensions of the rationals in which $2$ splits completely.

\section{Further results and questions}
\label{further-questions}

In this section we present some more computational observations
for Hecke algebras under the assumptions of Question~\ref{ourquestion},
which lead us to ask some more questions.

\subsection*{On the dimension of the Hecke algebra}

From the data, we see that many even integers appear as dimensions of
the~$\TT_\m$. We know that the dimension must be at least~4, as this
is the dimension of the smallest non-Gorenstein algebra which can
appear in our case. This extends the results
of~\cite{kilford-nongorenstein}, where the dimensions of the 
Hecke algebras~$\TT_{\ZZ \to \FF_2}(S_2(\Gamma_0(431)))$ 
and~$\TT_{\ZZ \to \FF_2}(S_2(\Gamma_0(503)))$ localised
at the non-Gorenstein maximal ideals are shown to be~4.

In this table we see exactly how many times each dimension appears in
our data. We observe that every even integer between~4 and~32 appears,
and that the largest dimension is~60. The most common dimension is~4,
which appears about half of the time. However, as the dimension of the
Hecke algebra attached to~$S_k(\Gamma_1(N))$ increases with~$N$ and
with~$k$, this may be an artifact of the data being collected for
``small'' levels~$N$ and primes~$p$.  \medskip

\begin{tabular}{|l|l|l|l|l|l|l|l|l|l|l|}
\hline
Dimension & 4  & 6  & 8  & 10 & 12 & 14 & 16 & 18 & 20 & 22  \\ %
\hline
Number of algebras &206 &  58& 25 &  3 & 24 &  6 & 20 &  3 & 12 &  3  \\ %
\hline
\hline
Dimension& 24 & 26 & 28 & 30 & 32 & 36 & 40 & 46 & 56 & 60\\
\hline
Number of algebras &5 &  4 &  2 &  1 & 2 &  2 &  4 &  1 &  2 &  1\\
\hline
\end{tabular}
\medskip

It seems reasonable that there should be infinitely many cases with
dimension~4, and plausible that every even integer greater than or
equal to~4 should appear as a dimension infinitely many times. From
the tables, we see that dimension~4 algebras appear at very high
levels, so they do not appear to be becoming rare as the dimension
increases, but this may, of course, be an artifact of our data.

We note that not every example that arises from an elliptic curve in
characteristic $p=2$ has Hecke algebra with dimension~$4$; for example
the algebra $\TT_{\ZZ \to \FF_2}(S_2(\Gamma_0(2089)))$ localised at
its non-Gorenstein maximal ideal has dimension~18. In level~$18097$
there is a dimension~$36$ example arising from an elliptic curve.

\subsection*{On the residue degree}

We will now solve an easy aspect of the question of the possible
structures of non-Gorenstein local algebras occurring as local Hecke
algebras. We assume for the couple of lines to follow the Generalised
Riemann Hypothesis (GRH).

We claim that then the residue degrees of $\TT_{\mathfrak{m}}$ (in the
notation of Question~\ref{ourquestion}) are unbounded, if we let~$p$ 
and $N$ run through the primes such that $p \neq N$ and $N$
is congruent to~$3$ modulo~$4$.

For, class groups of imaginary quadratic fields $\QQ(\sqrt{-N})$ have
arbitrarily large cyclic factors of odd order, as the exponent of
these class groups is known to go to infinity as $N$ does, by the main
result of~\cite{boyd}, which assumes GRH. So the discussion on
dihedral forms in Section~\ref{computational-results} immediately
implies the claim.

\subsection*{On the embedding dimension}

One can ask whether the embedding dimension of the local Hecke
algebras in the situation of Question~\ref{ourquestion} is bounded, if
we allow $p$ and~$N$ to vary. This, however, seems to be a difficult
problem. The embedding dimensions occuring in our tables are $3$ (299
times), $4$ (78 times) and~$5$ (7 times).

The embedding dimension~$d$ is related to the number of Hecke
operators needed to generate the local Hecke algebra, in the sense
that at least~$d$ Hecke operators are needed. Probably, $d$ Hecke
operators can be found that do generate, but they need not be the
first~$d$ prime Hecke operators, of course. However, as our tables
suggest, in most cases the actual computations were done using very few
operators, and there are 99 of the 384 cases when the computation
already finished after $d$ operators.

\section{Acknowledgements}

The authors would like to thank Kevin Buzzard for informing them about
the level 23 and weight 59 example, for many helpful conversations and
communications, and for his comments on a first draft of this paper.
They would like to thank Amod Agashe for comments on an early version
and William Stein for the use of his computers {\sc Meccah} and 
{\sc Neron} for computations; without these, they would not have been
able to pursue these computations as far as they did.  The authors are
grateful to Gunter Malle for providing them with a very useful list of
$A_5$-polynomials. The first author would also like to thank Edray
Goins for helpful conversations. The second author is endebted to Bas
Edixhoven for very enlightening explanations.

\bibliographystyle{plain}
\bibliography{nongorenstein}

\noindent \begin{tabular}{p{8cm}p{6cm}}
L.~J.~P.~Kilford          & Gabor Wiese             \\
Mathematics Institute     & NWF 1-Mathematik        \\
Oxford University         & Universität Regensburg  \\
24--29 St Giles'          & D-93040 Regensburg      \\
Oxford OX1 3LB            & Germany                 \\
United Kingdom            &                         \\
                          & {\tt gabor@pratum.net}  \\
{\tt l.kilford@gmail.com} & {\tt http://maths.pratum.net}
\end{tabular}

\appendix

\newpage

\section{The {\sc Magma} package HeckeAlgebra (by Gabor Wiese)}
\label{section-manual}

\abstract{This is a short manual for the {\sc Magma} package
{\tt HeckeAlgebra}, which can be downloaded from the author's
webpage.
The author would like to thank Lloyd Kilford for very helpful
suggestions.}

\subsection{Example}

{\setlength{\parindent}{0pt}

The following example explains the main functions of
the package.
Let us suppose that the file {\tt HeckeAlgebra.mg} is stored in 
the current path.
We first attach the package.\\
\magma{> Attach("HeckeAlgebra.mg");}

We want the package to be silent, so we put:\\
\magma{> SetVerbose ("HeckeAlgebra",false);}

If we would like more information on the computations being
performed, we should have put the value \magma{true}.
Since we want to store the data to be computed in a file,
we now create the file.

\magma{> my\_file := "datafile";\\
> CreateStorageFile(my\_file);}

Next, we would like to compute the Hecke algebras of the
dihedral eigenforms of level $2039$ over extensions of $\FF_2$.
First, we create a list of such forms.\\
\magma{> dih  := DihedralForms(2039 : ListOfPrimes := [2], completely\_split := 
false);}\\
Now, we compute the corresponding Hecke algebras, print part
of the computed data in a human readable format, and
finally save the data to our file.\\
\magma{> for f in dih do\\
for>     ha := HeckeAlgebras(f);\\
for>     HeckeAlgebraPrint1(ha);\\
for>     StoreData(my\_file, ha);\\
for> end for;}

\magmaout{
Level 2039\\
Weight 2\\
Characteristic 2\\
Gorenstein defect 0\\
Dimension 1\\
Number of operators used 3\\
Primes lt Hecke bound 68\\
Residue degree 2\\
---------------------------\\
Level 2039\\
Weight 2\\
Characteristic 2\\
Gorenstein defect 2\\
Dimension 6\\
Number of operators used 4\\
Primes lt Hecke bound 68\\
Residue degree 2\\
---------------------------\\
Level 2039\\
Weight 2\\
Characteristic 2\\
Gorenstein defect 0\\
Dimension 1\\
Number of operators used 3\\
Primes lt Hecke bound 68\\
Residue degree 6\\
---------------------------\\
Level 2039\\
Weight 2\\
Characteristic 2\\
Gorenstein defect 0\\
Dimension 1\\
Number of operators used 3\\
Primes lt Hecke bound 68\\
Residue degree 4\\
---------------------------\\
Level 2039\\
Weight 2\\
Characteristic 2\\
Gorenstein defect 0\\
Dimension 1\\
Number of operators used 3\\
Primes lt Hecke bound 68\\
Residue degree 4\\
---------------------------\\
Level 2039\\
Weight 2\\
Characteristic 2\\
Gorenstein defect 0\\
Dimension 1\\
Number of operators used 3\\
Primes lt Hecke bound 68\\
Residue degree 12\\
---------------------------\\
Level 2039\\
Weight 2\\
Characteristic 2\\
Gorenstein defect 0\\
Dimension 1\\
Number of operators used 3\\
Primes lt Hecke bound 68\\
Residue degree 12\\
---------------------------
}

With the function \magma{DihedralForms} one may also compute 
exclusively
representations that are completely split in the characteristic.
The default is \magma{completely\_split := true}. By the option
\magma{bound} we indicate primes up to which bound should be
used as the characteristic.
The following example illustrates this.\\
\magma{> dih1 := DihedralForms (431 : bound := 20);\\
> for f in dih1 do\\
for>    ha := HeckeAlgebras(f);\\
for>    HeckeAlgebraPrint1(ha);\\
for>    StoreData(my\_file, ha);\\
for> end for;}

\magmaout{
Level 431\\
Weight 2\\
Characteristic 2\\
Gorenstein defect 2\\
Dimension 4\\
Number of operators used 6\\
Primes lt Hecke bound 20\\
Residue degree 1\\
---------------------------\\
Level 431\\
Weight 11\\
Characteristic 11\\
Gorenstein defect 2\\
Dimension 4\\
Number of operators used 5\\
Primes lt Hecke bound 77\\
Residue degree 3\\
---------------------------
}

One can also compute icosahedral modular forms over extensions
of $\FF_2$, starting from an integer polynomial with
Galois group $A_5$, as follows.\\
\magma{> R<x> := PolynomialRing(Integers());\\
> pol := x\^{}5-x\^{}4-780*x\^{}3-1795*x\^{}2+3106*x+344;\\
> f := A5Form(pol);}

With this kind of icosahedral examples one has to pay attention to the
conductor, as it can be huge. This polynomial has prime conductor.
But conductors need not be square-free, in general.\\
\magma{> print Modulus(f`Character);}

\magmaout{1951}

So it's reasonable. We do the computation.\\
\magma{> ha := HeckeAlgebras(f);\\
> HeckeAlgebraPrint1(ha);}

\magmaout{
Level 1951\\
Weight 2\\
Characteristic 2\\
Gorenstein defect 0\\
Dimension 3\\
Number of operators used 3\\
Primes lt Hecke bound 66\\
Residue degree 4\\
---------------------------\\
Level 1951\\
Weight 2\\
Characteristic 2\\
Gorenstein defect 0\\
Dimension 6\\
Number of operators used 3\\
Primes lt Hecke bound 66\\
Residue degree 4\\
---------------------------
}

There are two forms, which is okay,
since they come from a weight one form in
two different ways and this case is not exceptional.
We now save them, as always.\\
\magma{> StoreData(my\_file, ha);}

It is also possible to compute all forms at a given character and weight.\\
\magma{> eps := DirichletGroup(229,GF(2)).1;\\
> ha := HeckeAlgebras(eps,2);\\
> HeckeAlgebraPrint1(ha);}

\magmaout{
Level 229\\
Weight 2\\
Characteristic 2\\
Gorenstein defect 0\\
Dimension 1\\
Number of operators used 12\\
Primes lt Hecke bound 12\\
Residue degree 1\\
---------------------------\\
Level 229\\
Weight 2\\
Characteristic 2\\
Gorenstein defect 0\\
Dimension 2\\
Number of operators used 12\\
Primes lt Hecke bound 12\\
Residue degree 2\\
---------------------------\\
Level 229\\
Weight 2\\
Characteristic 2\\
Gorenstein defect 0\\
Dimension 4\\
Number of operators used 12\\
Primes lt Hecke bound 12\\
Residue degree 1\\
---------------------------\\
Level 229\\
Weight 2\\
Characteristic 2\\
Gorenstein defect 0\\
Dimension 2\\
Number of operators used 12\\
Primes lt Hecke bound 12\\
Residue degree 5\\
---------------------------
}

\magma{> StoreData(my\_file,ha);}

Next, we illustrate how one reloads what has been saved.
One would like to type: \magma{load my\_file;}
but that does not work. One has to do it as follows.\\
\magma{> load "datafile";\\
> mf := RecoverData(LoadIn,LoadInRel);}

Now, \magma{mf} contains a list of all algebra data computed before.
There's a rather concise printing function, displaying part of the information,
namely \magma{HeckeAlgebraPrint(mf);}.

One can also create a LaTeX longtable. The entries can be chosen
in quite a flexible way. The standard usage is the following.\\
\magma{> HeckeAlgebraLaTeX(mf,"table.tex");}

A short LaTeX file displaying the table is the following:\\
{\tt {$\backslash$}documentclass[11pt]\{article\}\\
{$\backslash$}usepackage\{longtable\}\\
{$\backslash$}begin\{document\}\\
{$\backslash$}input\{table\}\\
{$\backslash$}end\{document\}}

The table we created is this one:
\begin{longtable}{||c|c|c|c|c|c|c|c|c|c||}
\hline
Level & Wt & ResD & Dim & EmbDim & NilO & GorDef & \#Ops & \#(p$<$HB) & Gp \\
\hline\endhead\hline\endfoot\hline\hline\endlastfoot
2039 & 2 & 2 & 1 & 0 & 0 & 0 & 3 & 68 & $D_{3}$ \\
2039 & 2 & 2 & 6 & 3 & 2 & 2 & 4 & 68 & $D_{5}$ \\
2039 & 2 & 6 & 1 & 0 & 0 & 0 & 3 & 68 & $D_{9}$ \\
2039 & 2 & 4 & 1 & 0 & 0 & 0 & 3 & 68 & $D_{15}$ \\
2039 & 2 & 4 & 1 & 0 & 0 & 0 & 3 & 68 & $D_{15}$ \\
2039 & 2 & 12 & 1 & 0 & 0 & 0 & 3 & 68 & $D_{45}$ \\
2039 & 2 & 12 & 1 & 0 & 0 & 0 & 3 & 68 & $D_{45}$ \\
431 & 2 & 1 & 4 & 3 & 1 & 2 & 6 & 20 & $D_{3}$ \\
431 & 11 & 3 & 4 & 3 & 1 & 2 & 5 & 77 & $D_{7}$ \\
1951 & 2 & 4 & 3 & 1 & 2 & 0 & 3 & 66 & $A_5$ \\
1951 & 2 & 4 & 6 & 2 & 3 & 0 & 3 & 66 & $A_5$ \\
229 & 2 & 1 & 1 & 0 & 0 & 0 & 12 & 12 & $$ \\
229 & 2 & 2 & 2 & 1 & 1 & 0 & 12 & 12 & $$ \\
229 & 2 & 1 & 4 & 1 & 3 & 0 & 12 & 12 & $$ \\
229 & 2 & 5 & 2 & 1 & 1 & 0 & 12 & 12 & $$ \\
\end{longtable}

In the examples of level $229$ the image of the Galois representation
as an abstract group is not know. That is due to the fact that we
created these examples without specifying the Galois representation in
advance.

It is possible to compute arbitrary Hecke operators
on the local Hecke factors generated by \magma{HeckeAlgebras($\cdot$)},
as the following example illustrates.\\
\magma{> A,B,M,C := HeckeAlgebras(DirichletGroup(253,GF(2)).1,2 : over\_residue\_field := true);}

Suppose that we want to know the Hecke operator $T_{17}$ on the $4$th local factor.\\
\magma{> i := 4;\\
> T := BaseChange(HeckeOperator(M,17),C[i]);}

The coefficients are the eigenvalues (only one):\\
\magma{> Eigenvalues(T);}

\magmaout{\{ $<$ \$ . 1 $\hat{ }$ 5, 8 $>$ \}}

Let us remember the eigenvalue.\\
\magma{> e := SetToSequence(Eigenvalues(T))[1][1];}

In order to illustrate the option \magma{over\_residue\_field}, 
we also compute the following:\\
\magma{> A1,B1,M1,C1 := HeckeAlgebras(DirichletGroup(253,GF(2)).1,2 : over\_residue\_field := false);\\
> T1 := BaseChange(HeckeOperator(M1,17),C1[i]);\\
> Eigenvalues(T1);}

\magmaout{\{\}}

The base field is strictly smaller than the residue field in this example
and the operator \magma{T1} cannot be diagonalised over the base field.
We check that \magma{e} is nevertheless a zero of the minimal 
polynomial of \magma{T1}.\\
\magma{> Evaluate(MinimalPolynomial(T1),e);}

\magmaout{0}

\bigskip

The precise usage of the package is described in the following
sections.
}

\subsection{Hecke algebra computation}

\subsubsection{The modular form format}\label{secmodforms}

In the package, modular forms are often represented by the following
record.\\
\noindent\magma{ModularFormFormat := recformat <}
\begin{longtable}{ll}
\magma{Character}           &\magma{: GrpDrchElt,}\\
\magma{Weight}              &\magma{: RngIntElt,}\\
\magma{CoefficientFunction} &\magma{: Map, }\\
\magma{ImageName}           &\magma{: MonStgElt,}\\ 
\magma{Polynomial}          &\magma{: RngUPolElt }
\end{longtable}
\magma{>;}

The fields \magma{Character} and \magma{Weight} have the obvious
meaning. Sometimes, the image of the associated Galois representation
is known as an abstract group. Then that name is recorded
in \magma{ImageName}, e.g.\ \magma{A\_5} or \magma{D\_3}.
In some cases, a polynomial is known whose splitting field
is the number field cut out by the Galois representation. Then
the polynomial is stored in \magma{Polynomial}. The cases in
which polynomials are known are usually icosahedral ones.
The \magma{CoefficientFunction} is a function 
from the integers to a polynomial ring. For all primes $l$ different
from the characteristic and not dividing the level of the
modular form (i.e.\ the modulus of the \magma{Character}), the 
coefficient function should return the minimal polynomial of
the $l$-th coefficient in the $q$-expansion of the modular
form in question.

\subsubsection{Dihedral modular forms}\label{secdih}

Eigenforms whose associated Galois representations take
dihedral groups as images provide an important source of examples,
in many contexts. These eigenforms are called {\em dihedral}.
The big advantage is that their Galois representation, and hence
their $q$-coefficients, can be computed using class field theory.
That enables one to exhibit Galois representations in the context
of modular forms with certain number theoretic properties.
The property for which these functions were initially created
is that the representations should be unramified in the 
characteristic, say $p$, and that $p$ is completely split
in the number field cut out by the representation.

We consider dihedral representations whose determinant
is the Legendre symbol of a quadratic field $\QQ(\sqrt{N})$. 
The representations produced by the functions to be described are obtained 
by induction of an unramified character $\chi$ of~$\QQ(\sqrt{N})$
whose conjugate by the non-trivial element of 
the Galois group of $\QQ(\sqrt{N})$ over $\QQ$ 
is assumed to be $\chi^{-1}$.

\intr{intrinsic GetLegendre (N :: RngIntElt, K :: FldFin ) -> GrpDrchElt}
For an odd positive integer \magma{N}, this function returns the element
of \magma{DirichletGroup(Abs(N),K)} (with \magma{K} a finite field of 
characteristic different from $2$) which
corresponds to the Legendre symbol $p \mapsto \left(\frac{\pm N}{p}\right)$.
If \magma{N} is $1$ mod $4$ the sign is $+1$, and $-1$ otherwise.

\intr{intrinsic DihedralForms (N :: RngIntElt : \\
\hspace*{2cm}        ListOfPrimes := [], 
                     bound := 100,              
                     odd\_only := true,
                     quad\_disc := 0, \\
\hspace*{2cm}        completely\_split := true,
                     all\_conjugacy\_classes := true
   ) -> Rec}
This function computes all modular forms (in the sense of
Section~\ref{secmodforms}) of level \magma{N} and
weight~$p$ over a 
finite field of characteristic $p$ that come from dihedral representations whose
determinant is the Legendre symbol of the quadratic field 
$K=\mathbb{Q}(\sqrt{\pm \text{\magma{quad\_disc}}})$ and which
are obtained by induction of an unramified character of~$K$.  
If \magma{quad\_disc} is $1$ mod $4$ the sign is $+1$, and $-1$ otherwise.
If \magma{quad\_disc} is $0$, the value of \magma{N} is used.
If the option \magma{completely\_split} is set, only those
representations are returned which are completely split at~$p$.  If
the option \magma{ListOfPrimes} is assigned a non-empty list of
primes, only those primes are considered as the characteristic.  If
it is the empty set, all primes $p$ up to the \magma{bound} are taken
into consideration.
If the option \magma{odd\_only} is true, only odd Galois representations
are returned.
If the option \magma{all\_conjugacy\_classes} is true, each unramified character 
as above up to Galois conjugacy and up to taking inverses is used. 
Otherwise, a single choice is made.
That there may be non-conjugate characters cutting out the same number
field is due to the fact that there may be non-conjugate elements
of the same order in the multiplicative group of a finite
field.

\subsubsection{Icosahedral modular forms}\label{secicosa}

Eigenforms whose attached Galois representations take the group $A_5$
as projective images are called {\em icosahedral}.  Since extensive
tables of $A_5$-extensions of the rationals are available, one can
consider icosahedral Galois representations which one knows very well.
That allows one to test certain conjectures concerning modular forms
on icosahedral ones.

We note the isomorphism $A_5 \cong \SL_2(\FF_4)$. Thus, $A_5$-extentions
of the rationals give rise to icosahedral Galois representations
in characteristic $2$ which (should) come from modular forms mod $2$.
It would also be possible to use certain other primes, but this has not
been implemented.

\intr{intrinsic A5Form (f :: RngUPolElt) -> Rec} 
Returns the icosahedral form in characteristic $2$ and weight $2$ of
smallest predicted level corresponding to the polynomial \magma{f}
which is expected to be of degree~$5$ and whose Galois group is
supposed to be~$A_5$.  No checks about \magma{f} are performed.

\subsubsection{The Hecke algebra format}

The data concerning the Hecke algebra of an eigenform
that is computed by the function \magma{HeckeAlgebras}
is a record of the following form.

\noindent\magma{AlgebraData := recformat <}
\begin{longtable}{ll}
  \magma{Level}              &\magma{: RngIntElt,}\\
  \magma{Weight}             &\magma{: RngIntElt,}\\
  \magma{Characteristic}     &\magma{: RngIntElt,}\\
  \magma{BaseFieldDegree}    &\magma{: RngIntElt,}\\
  \magma{CharacterOrder}     &\magma{: RngIntElt,}\\
  \magma{CharacterConductor} &\magma{: RngIntElt,}\\
  \magma{CharacterIndex}     &\magma{: RngIntElt,}\\
  \magma{AlgebraFieldDegree} &\magma{: RngIntElt,}\\
  \magma{ResidueDegree}      &\magma{: RngIntElt,}\\
  \magma{Dimension}          &\magma{: RngIntElt,}\\
  \magma{GorensteinDefect}   &\magma{: RngIntElt,}\\
  \magma{EmbeddingDimension} &\magma{: RngIntElt,}\\
  \magma{NilpotencyOrder}    &\magma{: RngIntElt,}\\
  \magma{Relations}          &\magma{: Tup,}\\
  \magma{NumberGenUsed}      &\magma{: RngIntElt,}\\
  \magma{ImageName}          &\magma{: MonStgElt,}\\
  \magma{Polynomial}         &\magma{: RngUPolElt}
\end{longtable}
\magma{>;}

\magma{Level} and \magma{Weight} have the obvious meaning.  Let $K$ be
the base field for the space of modular symbols used. It is (expected
to be) a finite field. Then \magma{Characteristic} is the
characteristic of $K$ and \magma{BaseFieldDegree} is the degree of $K$
over its prime field.  The entries \magma{CharacterOrder},
\magma{CharacterConductor} and \magma{CharacterIndex} concern the
Dirichlet character for which the modular symbols have been computed.
The latter field is the index of the character in
\magma{Elements(DirichletGroup($\cdot$))}. Note that that might change
between different versions of {\sc Magma}. The fields \magma{ResidueDegree}
(over the prime field), \magma{Dimension} and \magma{GorensteinDefect}
have their obvious meaning for the Hecke algebra in question.  The
tuple
\begin{quote}
\magma{<AlgebraFieldDegree, EmbeddingDimension, NilpotencyOrder,
Relations>}
\end{quote}
are data from which \magma{AffineAlgebra} can recreate the Hecke
algebra up to isomorphism.  \magma{NumberGenUsed} indicates the number
of generators used by the package for the computation of the Hecke
algebra. This number is usually much smaller than the Sturm bound.
\magma{ImageName} and \magma{Polynomial} have the same meaning as in
the record \magma{ModularFormFormat}.

\subsubsection{Hecke algebras}

\intr{intrinsic HeckeAlgebras (eps :: GrpDrchElt, weight :: RngIntElt :\\
\hspace*{2cm}      
    UserBound := 0,   
    first\_test := 3,      
    test\_interval := 1,   
    when\_test\_p := 3, \\    
\hspace*{2cm}      
    when\_test\_bad := 4,
    test\_sequence := [],
    dimension\_factor := 2,\\
\hspace*{2cm}
    ms\_space := 0,        
    cuspidal := true,
    DegreeBound := 0,
    OperatorList := [], \\
\hspace*{2cm}
    over\_residue\_field := true,
    try\_minimal := true,
    force\_local := false, \\  
\hspace*{.5cm}  ) -> SeqEnum, SeqEnum, ModSym, Tup, Tup\\
intrinsic HeckeAlgebras ( t :: Rec :\\
\hspace*{2cm}      
    UserBound := 0,   
    first\_test := 3,      
    test\_interval := 1,   
    when\_test\_p := 3, \\    
\hspace*{2cm}      
    when\_test\_bad := 4,
    test\_sequence := [],
    dimension\_factor := 2,\\
\hspace*{2cm}
    ms\_space := 0,        
    cuspidal := true,
    DegreeBound := 0,
    OperatorList := [], \\
\hspace*{2cm}
    over\_residue\_field := true,
    try\_minimal := true,
    force\_local := false, \\  
\hspace*{.5cm}  ) -> SeqEnum, SeqEnum, ModSym, Tup, Tup}

These functions compute all local Hecke algebras (up to Galois
conjugacy) in the specified \magma{weight} for the given Dirichlet
character \magma{eps}, respectively those corresponding to the modular
form \magma{t} given by a record of type \magma{ModularFormFormat}.
The functions return 5 values \magma{A,B,C,D,E}.  \magma{A} contains a
list of records of type \magma{AlgebraData} describing the local Hecke
algebra factors.  \magma{B} is a list containing the local Hecke
algebra factors as matrix algebras.  \magma{C} is the space of modular
symbols used in the computations.  \magma{D} is a tuple containing the
base change tuples describing the local Hecke factors. We need to know
\magma{D} in order to compute matrices representing Hecke operators
for the local factor.  Finally, \magma{E} contains a tuple consisting
of all Hecke operators computed so far for each local factor of the
Hecke algebra.

The usage in practice is described in the example at the beginning of
this manual. We now explain the different options in detail.

The modular symbols space to be used in the computations can be
determined as follows. The option \magma{ms\_space} can be set to the
values $1$ (the plus-space), $-1$ (the minus-space) and $0$ (the full
space).  Whether the restriction to the cuspidal subspace is taken, is
determined by \magma{cuspidal}. It is not necessary to pass to the
cuspidal subspace, for example, if a cusp form is given by a
coefficient function (see the description of the record
\magma{ModularFormFormat}).

In some cases, a list of Hecke operators on the modular symbols space
in question may already have been computed.  In order to prevent {\sc Magma}
from redoing their computations, they may be passed on to the function
using the option \magma{OperatorList}.

Often, one wants to compute the local Hecke algebra of a modular form
whose degree of the coefficient field over its prime field is known,
e.g.\ in the case of an icosahedral form in characteristic $2$ for the
trivial Dirichlet character the coefficient field is $\FF_4$.  By
assigning a positive value to the option \magma{DegreeBound} the
function will automatically discard any systems of eigenvalues beyond
that bound, which speeds up the computations.  One must be a bit
careful with this option, as there may be cases when the bound may not
be respected at ``bad primes''.  But it usually suffices to take twice
the degree of the coefficient field, e.g.\ one chooses
\magma{DegreeBound := 4} in the icosahedral example just described.
If no system of eigenvalues should be discarded for degree reasons,
one must set \magma{DegreeBound := 0}.

All of the options \magma{first\_test}, \magma{test\_interval},
\magma{when\_test\_p},\magma{when\_test\_bad}, \magma{test\_sequence},
\magma{force\_local}, \magma{dimension\_factor} and \magma{UserBound}
concern the stop criterion.  Theoretically, the Sturm bound (see
\magma{HeckeBound}) tells us up to which bound Hecke operators must be
computed in order to be sure that they generate the whole Hecke
algebra.  In practice, however, the algorithm can often determine
itself when enough Hecke operators have been computed to generate the
algebra. That number is usually much smaller than the Sturm bound. The
Sturm bound can be overwritten by assigning a positive number to
\magma{UserBound}.

The stop criterion is the following.  Let $M$ be the modular symbols
space used and $S$ the set of Hecke operators computed so far. Then $M
= \bigoplus_{i=1}^r M_i$ (for some $r$) such that each $M_i$ is
respected by the Hecke operators and the minimal polynomial of each $T
\in S$ restricted to $M_i$ is a power of an irreducible polynomial
(i.e.\ each $M_i$ is a primary space for the action of the algebra
generated by all elements of $S$). Let $A_i$ be the algebra generated
by $T|_{M_i}$ for all $T \in S$.  One knows (in many cases, and in all
cases of interest) that $A_i$ is equal to a direct product of local
Hecke algebras if one has the equality
$$ f \times \dim (A_i) = \text{ dimension of } M_i.$$
Here, $f$ is given by \magma{dimension\_factor} and should be $1$ if
the plus-space or the minus space of modular symbols are used, and $2$
otherwise.  The correct assignment of \magma{dimension\_factor} must
be made by hand, whence experimentations are possible.  If the stop
criterion is not reached, the algorithm terminates at the Hecke bound.

It may happen that, when the stop criterion is reached, one $A_i$ is
isomorphic to a direct product of more than one local Hecke algebras.
If in that case the option \magma{force\_local} is \magma{true}, the
computation of Hecke operators is continued until each $A_i$ is
isomorphic to a single Hecke factor. If \magma{force\_local} is
\magma{false}, then a fast localisation algorithm is applied to each
$A_i$.  The option is useful, when one expects only a single local
Hecke algebra factor, for example, when a modular form is given.

In many cases of interest the Hecke operator $T_p$ with $p$ the
characteristic is needed in order to generate the whole Hecke algebra.
The option \magma{when\_test\_p} tells the algorithm at which step to
compute $T_p$. It is very advisable to choose a small number. In
practice, the stop criterion is reached after very few steps, e.g.\ 5
steps, when $T_p$ is computed early. Otherwise, the algorithm often
has to continue until $T_p$ is computed, although most of the
operators before did not change the generated algebra. The option
\magma{when\_test\_bad} has a similar meaning for the $T_l$ for primes
$l$ dividing the level. However, paying attention to them is only
required when the modular form is old at~$l$. Moreover, one can assign
a list of primes to \magma{test\_sequence}. The algorithm will then
start with the Hecke operators indicated by that sequence, and then
continue with the others.

The option \magma{first\_test} tells the algorithm at which step the
first test for the stop criterion is to be performed. The next test is
then carried out after \magma{test\_interval} many steps, and so on.
These numbers should be chosen small, too, unless the dimension test
takes much time, which is rare, so that one wants to perform it less
often, meaning that possibly more Hecke operators than necessary are
computed (time consuming).

The option \magma{over\_residue\_field} tells the algorithm whether at
the end of the computation the local Hecke factors should be base
changed to their residue field. If that is done, only one of the
conjugate local factors of the base changed algebra is retained.

Finally, the option \magma{try\_minimal} is passed on to
\magma{AffineAlgebra}, when the output is generated. Calling that
function with the option set \magma{true} can sometimes be very time
consuming, but makes the output much shorter.

\subsubsection{Storage functions}

The package provides functions to store a list whose
elements are records of type \magma{AlgebraData} in a file, 
and to re-read it.
The usage of these functions is explained in the
example at the beginning of this manual.

\intr{intrinsic CreateStorageFile ( filename :: MonStgElt )}
This function prepares the file \magma{filename} for storing the data.

\intr{intrinsic StoreData (filename :: MonStgElt, forms :: SeqEnum)}
This functions appends the list \magma{forms} of Hecke algebra data 
to the file \magma{filename}. That
file must have been created by \magma{CreateStorageFile}.

\intr{intrinsic StoreData (filename :: MonStgElt, form :: Rec)}
This function appends the Hecke algebra data \magma{form}
to the file \magma{filename}. That
file must have been created by \magma{CreateStorageFile}.

\intr{intrinsic RecoverData (LoadIn :: SeqEnum, LoadInRel :: Tup ) -> SeqEnum}
In order to read Hecke algebra data from file \magma{``name''},
proceed as follows:\\
\magma{\hspace*{.5cm}> load ``name'';\\
\hspace*{.5cm}> readData := RecoverData(LoadIn,LoadInRel).}\\
Then \magma{readData} will contain a list whose elements are
records of type \magma{AlgebraData}.

\subsubsection{Output functions}

\intr{intrinsic HeckeAlgebraPrint (ha :: SeqEnum)\\
intrinsic HeckeAlgebraPrint1 (ha :: SeqEnum)}
These functions print part of the data stored in the list \magma{ha}
of records of type \magma{AlgebraData} in a human readable format.

\intr{intrinsic GetLevel (a :: Rec) -> Any\\
intrinsic GetWeight (a :: Rec) -> Any\\
intrinsic GetCharacteristic (a :: Rec) -> Any\\
intrinsic GetResidueDegree (a :: Rec) -> Any\\
intrinsic GetDimension (a :: Rec) -> Any\\
intrinsic GetGorensteinDefect (a :: Rec) -> Any\\
intrinsic GetEmbeddingDimension (a :: Rec) -> Any\\
intrinsic GetNilpotencyOrder (a :: Rec) -> Any\\
intrinsic GetHeckeBound (a :: Rec) -> Any\\
intrinsic GetPrimesUpToHeckeBound (a :: Rec) -> Any\\
intrinsic GetNumberOperatorsUsed (a :: Rec) -> Any\\
intrinsic GetPolynomial (a :: Rec) -> Any\\
intrinsic GetImageName (a :: Rec) -> Any}
These functions return the property of the record \magma{a}
of type \magma{AlgebraData} specified by the name of the
function. If the corresponding attribute is not assigned,
the empty string is returned.

\intr{intrinsic HeckeAlgebraLaTeX (ha :: SeqEnum, filename :: MonStgElt : which := [ \\
\hspace*{2cm}
             <GetLevel,"Level">,
             <GetWeight,"Wt">,
             <GetResidueDegree,"ResD">,\\
\hspace*{2cm}
             <GetDimension,"Dim">,
             <GetEmbeddingDimension,"EmbDim">,\\
\hspace*{2cm}
             <GetNilpotencyOrder,"NilO">,
             <GetGorensteinDefect,"GorDef">,\\
\hspace*{2cm}
             <GetNumberOperatorsUsed,"\#Ops">,\\
\hspace*{2cm}
             <GetPrimesUpToHeckeBound,"\#(p$<$HB)">,
             <GetImageName,"Gp">
                                                    ] )}
This function creates the LaTeX file \magma{filename} containing a
longtable consisting of certain properties of the objects in
\magma{ha} which are supposed to be records of type
\magma{AlgebraData}.  The properties to be written are indicated by
the list given in the option \magma{which} consisting of tuples
\magma{<f, name>}.  Here \magma{f} is a function that evaluates a
record of type \magma{AlgebraData} to some Magma object which is
afterwards transformed into a string using \magma{Sprint}. Examples
for \magma{f} are the functions \magma{GetLevel} etc., which are
described above.  The \magma{name} will appear in the table header.
For a sample usage, see the example at the beginning of this manual.

\subsubsection{Other functions}

\intr{intrinsic HeckeBound ( N :: RngIntElt, k :: RngIntElt ) -> RngIntElt\\
intrinsic HeckeBound ( eps :: GrpDrchElt, k :: RngIntElt ) -> RngIntElt}
These functions compute the Hecke bound for weight \magma{k} and
level \magma{N}, respectively Dirichelt character \magma{eps}.
Note that the Hecke bound is also often called the Sturm bound.

\subsection{Algebra handling}

\subsubsection{Affine algebras}

Let $A$ be a commutative local Artin algebra with maximal ideal $\fm$
over a finite field $k$. The residue field $K=A/\fm$ is a finite
extension of~$k$. By base changing to $K$ and taking one of the
conjugate local factors, we now assume that $k=K$.  The {\em embedding
  dimension} $e$ is the $k$-dimension of $\fm/\fm^2$.  By Nakayama's
Lemma, this is the minimal number of generators for $\fm$.  The name
comes from the fact that there is a surjection
$$    \pi: k[x_1, \dots, x_e] \twoheadrightarrow A.$$
Its kernel is called the {\em relations ideal}. By the {\em nilpotency
order} we mean the maximal integer $n$ such that $m^n$ is not the
zero ideal.  (As the algebra is local and Artin, its maximal ideal is
nilpotent.)  We know that the ideal
$$   J^{n+1} \text{ with } J := (x_1, \dots, x_e) $$
is in the kernel of $\pi$. 
So, in order to store $\pi$, we only need to store 
the kernel $R$ of the linear map between two finite dimensional $k$-vector 
spaces
$$   \pi_1: k[x_1, \dots, x_e]/J^{n+1} \twoheadrightarrow A.$$

From the tuple $<k,e,n,R>$ the algebra can be recreated (up to
isomorphism).  Let us point out, however, that from the tuple it is
not obvious whether two algebras are isomorphic. That would have to be
tested after recreating the algebras.

These functions are used in order to store the Hecke algebras computed
by \magma{HeckeAlgebras} in a way that does not use much memory, but
retains the algebra up to isomorphism.

\intr{intrinsic AffineAlgebra (A :: AlgMat : try\_minimal := true) -> RngMPolRes\\
intrinsic AffineAlgebra (A :: AlgAss : try\_minimal := true) -> RngMPolRes}
This function turns the local commutative algebra \magma{A} into an affine algebra
over its residue field. In fact, the algebra is first base changed
to its residue field, then for one of the conjugate local factors
an affine presentation is computed.
If the option \magma{try\_minimal} is true, the number of relations
will in general be smaller, but the computation time may be longer.

\intr{intrinsic AffineAlgebraTup (A :: AlgMat : try\_minimal := true ) -> Tup\\
intrinsic AffineAlgebraTup (A :: AlgAss : try\_minimal := true) -> Tup}
Given a commutative local Artin algebra \magma{A}, this function returns
a tuple \magma{<k,e,n,R>}, consisting of the
residue field \magma{k} of \magma{A}, the embedding dimension \magma{e}, 
the nilpotency order \magma{n} and relations \magma{R}.
From these data, an affine algebra can be recreated which is
isomorphic to one of the local factors of $A$ base changed
to its residue field.
If the option \magma{try\_minimal} is true, the number of relations
will in general be smaller, but the computation time may be longer.

\intr{intrinsic AffineAlgebra (form :: Rec) -> RngMPolRes}
Given a record of type \magma{AlgebraData}, this function returns the
corresponding Hecke algebra as an affine algebra.

\intr{intrinsic AffineAlgebra (A :: Tup) -> RngMPolRes}
This function turns a tuple \magma{<k,e,n,R>}, as above
consisting of a field \magma{k}, two integers
\magma{e}, \magma{n} (the embedding dimension and the nilpotency order) 
and relations \magma{R}, into an affine algebra.

\subsubsection{Matrix algebra functions}

\intr{intrinsic MatrixAlgebra ( L :: SeqEnum ) -> AlgMat}
Given a list of matrices \magma{L}, this function returns the matrix
algebra generated by the members of~\magma{L}.

\intr{intrinsic RegularRepresentation ( A :: AlgMat ) -> AlgMat}
This function computes the regular representation 
of the commutative matrix algebra \magma{A}.

\intr{intrinsic CommonLowerTriangular ( A :: AlgMat ) -> AlgMat}
Given a local commutative matrix algebra \magma{A}, this function
returns an isomorphic matrix algebra whose matrices are all lower
triangular, after a scalar extension to the residue field and taking
one of the Galois conjugate factors.

\bigskip
\noindent{\bf Base change}

\intr{intrinsic BaseChange ( S :: Tup, T :: Tup ) -> Tup}
This function computes the composition of the base change matrices
\magma{T = <C,D>}, followed by those in \magma{S = <E,F>}.

\intr{intrinsic BaseChange ( M :: Mtrx, T :: Tup ) -> Mtrx}
Given a matrix \magma{M} and a tuple \magma{T = <C,D>} of base change
matrices (for a subspace), this function computes the matrix of
\magma{M} with respect to the basis corresponding to \magma{T}.

\intr{intrinsic BaseChange ( M :: AlgMat, T :: Tup ) -> AlgMat}
Given a matrix algebra \magma{M} and a tuple \magma{T = <C,D>} of base
change matrices (for a subspace), this function computes the matrix
algebra of \magma{M} with respect to the basis corresponding to
\magma{T}.

\bigskip
\noindent{\bf Decomposition}

\intr{intrinsic Decomposition ( M :: Mtrx : DegBound := 0 ) -> Tup\\
intrinsic DecompositionUpToConjugation ( M :: Mtrx : DegBound := 0 ) -> Tup}
Given a matrix \magma{M}, these functions compute
a decomposition of the standard vector space such that 
\magma{M} acts as multiplication by a scalar 
on each summand. The output is a tuple consisting of
base change tuples \magma{<C,D>} corresponding to the summands.
With the second usage, summands conjugate under the absolute Galois group 
only appear once.

\intr{intrinsic Decomposition ( L :: SeqEnum : DegBound := 0 ) -> Tup\\
intrinsic DecompositionUpToConjugation ( L :: SeqEnum : DegBound := 0) -> Tup}
Given a sequence \magma{L} of commuting matrices, these
functions compute a decomposition of the standard vector space such that 
each matrix in \magma{L} acts as multiplication by a scalar 
on each summand. The output is a tuple consisting of
base change tuples \magma{<C,D>} corresponding to the summands.
With the second usage, summands conjugate under the absolute Galois group 
only appear once.

\intr{intrinsic Decomposition ( A :: AlgMat : DegBound := 0 ) -> Tup\\
intrinsic DecompositionUpToConjugation ( A :: AlgMat : DegBound := 0 ) -> Tup}
Given a commutative matrix algebra \magma{A}, these functions 
compute a decomposition of the standard vector space such that 
each element in \magma{A} acts as multiplication by a scalar 
on each summand. The output is a tuple consisting of
base change tuples \magma{<C,D>} corresponding to the summands.
With the second usage, summands conjugate under the absolute Galois group 
only appear once.

\intr{intrinsic AlgebraDecomposition ( A :: AlgMat : DegBound := 0 ) -> SeqEnum\\
intrinsic AlgebraDecompositionUpToConjugation ( A :: AlgMat : DegBound := 0 )\\
\hspace*{2cm} -> SeqEnum}
Given a matrix algebra \magma{A} over a finite field,
these functions return a local factor of \magma{A} after
scalar extension to the residue field.
With the second usage, factors conjugate under the absolute Galois group 
only appear once.

\intr{intrinsic ChangeToResidueField ( A :: AlgMat ) -> SeqEnum}
This function is identical to \magma{AlgebraDecompositionUpToConjugation}.

\bigskip
\noindent{\bf Localisations}

\intr{intrinsic Localisations ( L :: SeqEnum ) -> Tup, Tup\\
intrinsic Localisations ( A :: AlgMat ) -> Tup, Tup}
Given a list \magma{L} of commuting matrices or a commutative matrix
algebra \magma{A}, this function computes two tuples \magma{C},
\magma{D}, where \magma{C} contains a tuple consisting of the
localisations of \magma{A}, respectively of the matrix algebra
generated by \magma{L}, and \magma{D} consists of the corresponding
base change tuples.

\subsubsection{Associative algebras}

\intr{intrinsic Localisations ( A :: AlgAss ) -> SeqEnum}
This function returns a list of all localisations of the Artin algebra
\magma{A}, which is assumed to be commutative.  The output is a list
of associative algebras.

\subsubsection{Gorenstein defect}

Let $A$ be a local Artin algebra over a field with unique maximal
ideal $\fm$.  We define the {\em Gorenstein defect} of $A$ to be
$(\dim_{A/\fm} A[\fm]) -1$, which is equal to the number of $A$-module
generators of the annihilator of the maximal ideal minus one.  The
algebra is said to be {\em Gorenstein} if its Gorenstein defect is
equal to~$0$.

\intr{intrinsic GorensteinDefect ( A :: RngMPolRes) -> RngIntElt\\
intrinsic GorensteinDefect ( A :: AlgAss) -> RngIntElt\\
intrinsic GorensteinDefect ( A :: AlgMat ) -> RngIntElt}
These functions return the Gorenstein defect of the 
local commutative algebra \magma{A}.

\intr{intrinsic IsGorenstein ( M :: RngMPolRes ) -> BoolElt\\
intrinsic IsGorenstein ( M :: AlgAss ) -> BoolElt\\
intrinsic IsGorenstein ( M :: AlgMat ) -> BoolElt}
These functions test whether the commutative local algebra 
\magma{M} is Gorenstein.

\newpage

\section{Tables of Hecke algebras}
\label{section-tables}

\subsection*{Characteristic $p=2$, prime levels, dihedral}

\begin{longtable}{||c|c|c|c|c|c|c|c|c|c||}
\hline
Level & Wt & ResD & Dim & EmbDim & NilO & GorDef & \#Ops & \#(p$<$HB) & Gp \\
\hline\endhead\hline\endfoot\hline\hline\endlastfoot
431 {}\footnote{See \cite{kilford-nongorenstein}.} & 2 & 1 & 4 & 3 & 1 & 2 & 6 & 20 & $D_{3}$ \\
503 {}\footnote{See \cite{kilford-nongorenstein}.} & 2 & 1 & 4 & 3 & 1 & 2 & 3 & 23 & $D_{3}$ \\
1319 & 2 & 2 & 4 & 3 & 1 & 2 & 6 & 47 & $D_{5}$ \\
1439 & 2 & 1 & 4 & 3 & 1 & 2 & 4 & 52 & $D_{3}$ \\
1559 & 2 & 1 & 4 & 3 & 1 & 2 & 7 & 55 & $D_{3}$ \\
1607 & 2 & 1 & 4 & 3 & 1 & 2 & 3 & 56 & $D_{3}$ \\
1759 & 2 & 1 & 4 & 3 & 1 & 2 & 5 & 62 & $D_{3}$ \\
1823 & 2 & 2 & 4 & 3 & 1 & 2 & 3 & 62 & $D_{5}$ \\
1879 & 2 & 1 & 16 & 4 & 5 & 2 & 6 & 65 & $D_{3}$ \\
1951 & 2 & 1 & 4 & 3 & 1 & 2 & 4 & 66 & $D_{3}$ \\
1999 & 2 & 1 & 4 & 3 & 1 & 2 & 5 & 67 & $D_{3}$ \\
2039 & 2 & 2 & 6 & 3 & 2 & 2 & 4 & 68 & $D_{5}$ \\
2089 {}\footnote{See \cite{kilford-nongorenstein}.} & 2 & 1 & 18 & 4 & 7 & 2 & 5 & 70 & $D_{3}$ \\
2351 & 2 & 1 & 6 & 3 & 2 & 2 & 6 & 77 & $D_{3}$ \\
3407 & 2 & 1 & 16 & 4 & 5 & 2 & 7 & 103 & $D_{3}$ \\
3527 & 2 & 2 & 4 & 3 & 1 & 2 & 3 & 107 & $D_{5}$ \\
3623 & 2 & 1 & 6 & 3 & 2 & 2 & 3 & 110 & $D_{3}$ \\
3967 & 2 & 1 & 14 & 4 & 4 & 2 & 4 & 121 & $D_{3}$ \\
4231 & 2 & 1 & 4 & 3 & 1 & 2 & 10 & 126 & $D_{3}$ \\
4481 & 2 & 1 & 8 & 4 & 2 & 2 & 7 & 132 & $D_{3}$ \\
4799 & 2 & 1 & 12 & 4 & 3 & 2 & 5 & 139 & $D_{3}$ \\
4943 & 2 & 2 & 6 & 3 & 2 & 2 & 4 & 143 & $D_{5}$ \\
5167 & 2 & 1 & 6 & 3 & 2 & 2 & 5 & 149 & $D_{3}$ \\
5591 & 2 & 1 & 12 & 4 & 3 & 2 & 8 & 158 & $D_{3}$ \\
5591 & 2 & 3 & 4 & 3 & 1 & 2 & 5 & 158 & $D_{9}$ \\
5791 & 2 & 1 & 8 & 3 & 3 & 2 & 8 & 162 & $D_{3}$ \\
6199 & 2 & 1 & 16 & 4 & 5 & 2 & 7 & 174 & $D_{3}$ \\
6287 & 2 & 1 & 6 & 3 & 2 & 2 & 4 & 175 & $D_{3}$ \\
6343 & 2 & 1 & 12 & 4 & 3 & 2 & 5 & 177 & $D_{3}$ \\
6551 & 2 & 1 & 6 & 3 & 2 & 2 & 7 & 182 & $D_{3}$ \\
6823 & 2 & 1 & 4 & 3 & 1 & 2 & 4 & 189 & $D_{3}$ \\
6911 & 2 & 1 & 4 & 3 & 1 & 2 & 4 & 190 & $D_{3}$ \\
6967 & 2 & 1 & 12 & 4 & 3 & 2 & 8 & 191 & $D_{3}$ \\
7057 & 2 & 1 & 16 & 4 & 4 & 2 & 6 & 193 & $D_{3}$ \\
7103 & 2 & 3 & 4 & 3 & 1 & 2 & 3 & 194 & $D_{7}$ \\
7151 & 2 & 2 & 4 & 3 & 1 & 2 & 4 & 195 & $D_{5}$ \\
7351 & 2 & 1 & 12 & 4 & 3 & 2 & 9 & 200 & $D_{3}$ \\
7487 & 2 & 2 & 6 & 3 & 2 & 2 & 3 & 203 & $D_{5}$ \\
7583 & 2 & 1 & 4 & 3 & 1 & 2 & 6 & 205 & $D_{3}$ \\
7951 & 2 & 1 & 4 & 3 & 1 & 2 & 6 & 216 & $D_{3}$ \\
8111 & 2 & 5 & 4 & 3 & 1 & 2 & 5 & 217 & $D_{11}$ \\
8167 & 2 & 1 & 6 & 3 & 2 & 2 & 4 & 218 & $D_{3}$ \\
8191 & 2 & 2 & 6 & 3 & 2 & 2 & 5 & 218 & $D_{5}$ \\
8623 & 2 & 1 & 4 & 3 & 1 & 2 & 8 & 227 & $D_{3}$ \\
8713 & 2 & 1 & 16 & 4 & 4 & 2 & 4 & 231 & $D_{3}$ \\
9127 & 2 & 1 & 16 & 4 & 5 & 2 & 4 & 240 & $D_{3}$ \\
9281 & 2 & 1 & 12 & 4 & 3 & 2 & 8 & 243 & $D_{3}$ \\
9439 & 2 & 1 & 8 & 3 & 3 & 2 & 4 & 248 & $D_{3}$ \\
9623 & 2 & 2 & 6 & 3 & 2 & 2 & 5 & 252 & $D_{5}$ \\
9967 & 2 & 1 & 6 & 3 & 2 & 2 & 4 & 260 & $D_{3}$ \\
10079 & 2 & 1 & 12 & 4 & 3 & 2 & 8 & 263 & $D_{3}$ \\
10103 & 2 & 1 & 4 & 3 & 1 & 2 & 4 & 263 & $D_{3}$ \\
10391 & 2 & 1 & 6 & 3 & 2 & 2 & 9 & 269 & $D_{3}$ \\
10391 & 2 & 3 & 4 & 3 & 1 & 2 & 8 & 269 & $D_{9}$ \\
10487 & 2 & 1 & 10 & 3 & 4 & 2 & 3 & 272 & $D_{3}$ \\
10567 & 2 & 1 & 12 & 3 & 5 & 2 & 5 & 274 & $D_{3}$ \\
10639 & 2 & 1 & 4 & 3 & 1 & 2 & 4 & 274 & $D_{3}$ \\
10663 & 2 & 1 & 6 & 3 & 2 & 2 & 9 & 275 & $D_{3}$ \\
10687 & 2 & 1 & 6 & 3 & 2 & 2 & 4 & 275 & $D_{3}$ \\
10799 & 2 & 1 & 14 & 3 & 6 & 2 & 8 & 278 & $D_{3}$ \\
11159 & 2 & 3 & 4 & 3 & 1 & 2 & 4 & 283 & $D_{7}$ \\
11321 & 2 & 1 & 8 & 4 & 2 & 2 & 9 & 289 & $D_{3}$ \\
11743 & 2 & 1 & 4 & 3 & 1 & 2 & 5 & 297 & $D_{3}$ \\
13063 & 2 & 1 & 6 & 3 & 2 & 2 & 5 & 326 & $D_{3}$ \\
13487 & 2 & 1 & 8 & 3 & 3 & 2 & 6 & 334 & $D_{3}$ \\
13999 & 2 & 3 & 4 & 3 & 1 & 2 & 4 & 345 & $D_{7}$ \\
14303 & 2 & 1 & 4 & 3 & 1 & 2 & 5 & 354 & $D_{3}$ \\
14543 & 2 & 3 & 4 & 3 & 1 & 2 & 3 & 360 & $D_{7}$ \\
14639 & 2 & 2 & 4 & 3 & 1 & 2 & 5 & 361 & $D_{5}$ \\
15121 & 2 & 2 & 8 & 4 & 2 & 2 & 6 & 369 & $D_{5}$ \\
15193 & 2 & 1 & 16 & 4 & 6 & 2 & 4 & 370 & $D_{3}$ \\
15271 & 2 & 1 & 6 & 3 & 2 & 2 & 8 & 372 & $D_{3}$ \\
15383 & 2 & 2 & 6 & 3 & 2 & 2 & 3 & 375 & $D_{5}$ \\
15391 & 2 & 1 & 6 & 3 & 2 & 2 & 7 & 375 & $D_{3}$ \\
15551 & 2 & 1 & 4 & 3 & 1 & 2 & 11 & 377 & $D_{3}$ \\
15607 & 2 & 1 & 6 & 3 & 2 & 2 & 8 & 378 & $D_{3}$ \\
15641 & 2 & 1 & 32 & 4 & 7 & 2 & 8 & 378 & $D_{3}$ \\
15919 & 2 & 1 & 26 & 4 & 7 & 2 & 6 & 383 & $D_{3}$ \\
15991 & 2 & 1 & 12 & 4 & 3 & 2 & 9 & 386 & $D_{3}$ \\
16127 & 2 & 3 & 4 & 3 & 1 & 2 & 3 & 390 & $D_{7}$ \\
16369 & 2 & 1 & 24 & 4 & 6 & 2 & 6 & 398 & $D_{3}$ \\
16487 & 2 & 1 & 6 & 3 & 2 & 2 & 4 & 400 & $D_{3}$ \\
16649 & 2 & 1 & 16 & 4 & 4 & 2 & 8 & 403 & $D_{3}$ \\
17471 & 2 & 1 & 6 & 3 & 2 & 2 & 11 & 421 & $D_{3}$ \\
18047 & 2 & 1 & 30 & 4 & 6 & 2 & 4 & 431 & $D_{3}$ \\
18097 & 2 & 1 & 36 & 5 & 6 & 2 & 9 & 432 & $D_{3}$ \\
18127 & 2 & 1 & 12 & 4 & 3 & 2 & 5 & 433 & $D_{3}$ \\
18257 & 2 & 1 & 12 & 4 & 4 & 2 & 4 & 436 & $D_{3}$ \\
19079 & 2 & 1 & 4 & 3 & 1 & 2 & 4 & 449 & $D_{3}$ \\
19079 & 2 & 3 & 4 & 3 & 1 & 2 & 4 & 449 & $D_{9}$ \\
19441 & 2 & 1 & 16 & 4 & 4 & 2 & 7 & 457 & $D_{3}$ \\
19543 & 2 & 2 & 4 & 3 & 1 & 2 & 4 & 460 & $D_{5}$ \\
19583 & 2 & 1 & 4 & 3 & 1 & 2 & 7 & 461 & $D_{3}$ \\
19751 & 2 & 1 & 4 & 3 & 1 & 2 & 4 & 462 & $D_{3}$ \\
19919 & 2 & 1 & 16 & 4 & 5 & 2 & 7 & 467 & $D_{3}$ \\
19927 & 2 & 1 & 6 & 3 & 2 & 2 & 4 & 467 & $D_{3}$ \\
20183 & 2 & 2 & 4 & 3 & 1 & 2 & 7 & 474 & $D_{5}$ \\
20599 & 2 & 1 & 6 & 3 & 2 & 2 & 6 & 481 & $D_{3}$ \\
20759 & 2 & 1 & 18 & 4 & 6 & 2 & 9 & 483 & $D_{3}$ \\
20887 & 2 & 2 & 6 & 3 & 2 & 2 & 4 & 487 & $D_{5}$ \\
21319 & 2 & 3 & 4 & 3 & 1 & 2 & 9 & 497 & $D_{7}$ \\
21647 & 2 & 1 & 6 & 3 & 2 & 2 & 6 & 504 & $D_{3}$ \\
21737 & 2 & 1 & 24 & 5 & 4 & 2 & 7 & 507 & $D_{3}$ \\
21839 & 2 & 2 & 4 & 3 & 1 & 2 & 7 & 509 & $D_{5}$ \\
22159 & 2 & 2 & 4 & 3 & 1 & 2 & 9 & 515 & $D_{5}$ \\
22511 & 2 & 1 & 4 & 3 & 1 & 2 & 6 & 522 & $D_{3}$ \\
22567 & 2 & 3 & 4 & 3 & 1 & 2 & 4 & 523 & $D_{7}$ \\
22751 & 2 & 3 & 4 & 3 & 1 & 2 & 4 & 526 & $D_{7}$ \\
23159 & 2 & 1 & 20 & 4 & 6 & 2 & 6 & 535 & $D_{3}$ \\
23159 & 2 & 3 & 6 & 3 & 2 & 2 & 5 & 535 & $D_{9}$ \\
23279 & 2 & 1 & 8 & 3 & 3 & 2 & 4 & 537 & $D_{3}$ \\
23321 & 2 & 1 & 12 & 4 & 4 & 2 & 7 & 538 & $D_{3}$ \\
23417 & 2 & 1 & 26 & 4 & 7 & 2 & 10 & 539 & $D_{3}$ \\
23567 & 2 & 1 & 4 & 3 & 1 & 2 & 3 & 544 & $D_{3}$ \\
23687 & 2 & 3 & 4 & 3 & 1 & 2 & 3 & 548 & $D_{7}$ \\
23743 & 2 & 2 & 4 & 3 & 1 & 2 & 4 & 548 & $D_{5}$ \\
24151 & 2 & 2 & 4 & 3 & 1 & 2 & 5 & 556 & $D_{5}$ \\
24281 & 2 & 1 & 16 & 4 & 4 & 2 & 5 & 557 & $D_{3}$ \\
24439 & 2 & 2 & 6 & 3 & 2 & 2 & 7 & 561 & $D_{5}$ \\
24847 & 2 & 2 & 4 & 3 & 1 & 2 & 8 & 570 & $D_{5}$ \\
25031 & 2 & 1 & 8 & 3 & 3 & 2 & 7 & 573 & $D_{3}$ \\
25111 & 2 & 1 & 6 & 3 & 2 & 2 & 8 & 574 & $D_{3}$ \\
25247 & 2 & 1 & 6 & 3 & 2 & 2 & 10 & 575 & $D_{3}$ \\
25409 & 2 & 1 & 8 & 4 & 2 & 2 & 8 & 580 & $D_{3}$ \\
25439 & 2 & 1 & 4 & 3 & 1 & 2 & 6 & 580 & $D_{3}$ \\
25447 & 2 & 1 & 6 & 3 & 2 & 2 & 11 & 581 & $D_{3}$ \\
25793 & 2 & 1 & 16 & 4 & 4 & 2 & 5 & 590 & $D_{3}$ \\
26431 & 2 & 2 & 6 & 3 & 2 & 2 & 4 & 599 & $D_{5}$ \\
26839 & 2 & 3 & 4 & 3 & 1 & 2 & 5 & 607 & $D_{7}$ \\
26959 & 2 & 1 & 6 & 3 & 2 & 2 & 23 & 610 & $D_{3}$ \\
27143 & 2 & 1 & 4 & 3 & 1 & 2 & 4 & 615 & $D_{3}$ \\
27143 & 2 & 3 & 4 & 3 & 1 & 2 & 3 & 615 & $D_{9}$ \\
27647 & 2 & 1 & 26 & 4 & 10 & 2 & 7 & 623 & $D_{3}$ \\
27673 & 2 & 1 & 8 & 4 & 2 & 2 & 7 & 623 & $D_{3}$ \\
27743 & 2 & 3 & 4 & 3 & 1 & 2 & 3 & 624 & $D_{7}$ \\
28031 & 2 & 2 & 6 & 3 & 2 & 2 & 5 & 631 & $D_{5}$ \\
28031 & 2 & 1 & 8 & 3 & 3 & 2 & 5 & 631 & $D_{3}$ \\
28031 & 2 & 4 & 4 & 3 & 1 & 2 & 7 & 631 & $D_{15}$ \\
28279 & 2 & 1 & 8 & 3 & 3 & 2 & 4 & 635 & $D_{3}$ \\
28279 & 2 & 1 & 26 & 4 & 8 & 2 & 14 & 635 & $D_{3}$ \\
28279 & 2 & 1 & 4 & 3 & 1 & 2 & 11 & 635 & $D_{3}$ \\
28279 & 2 & 1 & 6 & 3 & 2 & 2 & 6 & 635 & $D_{3}$ \\
28703 & 2 & 1 & 20 & 4 & 5 & 2 & 5 & 642 & $D_{3}$ \\
28759 & 2 & 1 & 8 & 3 & 3 & 2 & 5 & 645 & $D_{3}$ \\
29023 & 2 & 1 & 4 & 3 & 1 & 2 & 8 & 650 & $D_{3}$ \\
29287 & 2 & 2 & 8 & 3 & 3 & 2 & 4 & 653 & $D_{5}$ \\
29311 & 2 & 3 & 6 & 3 & 2 & 2 & 5 & 653 & $D_{7}$ \\
29399 & 2 & 1 & 6 & 3 & 2 & 2 & 7 & 654 & $D_{3}$ \\
29567 & 2 & 1 & 6 & 3 & 2 & 2 & 6 & 657 & $D_{3}$ \\
29879 & 2 & 2 & 4 & 3 & 1 & 2 & 5 & 666 & $D_{5}$ \\
29959 & 2 & 1 & 6 & 3 & 2 & 2 & 4 & 668 & $D_{3}$ \\
29959 & 2 & 3 & 4 & 3 & 1 & 2 & 4 & 668 & $D_{9}$ \\
30223 & 2 & 1 & 8 & 3 & 3 & 2 & 5 & 674 & $D_{3}$ \\
30367 & 2 & 2 & 6 & 3 & 2 & 2 & 5 & 677 & $D_{5}$ \\
30431 & 2 & 2 & 4 & 3 & 1 & 2 & 4 & 677 & $D_{5}$ \\
30559 & 2 & 3 & 4 & 3 & 1 & 2 & 9 & 680 & $D_{7}$ \\
30727 & 2 & 1 & 6 & 3 & 2 & 2 & 9 & 685 & $D_{3}$ \\
30911 & 2 & 2 & 4 & 3 & 1 & 2 & 5 & 686 & $D_{5}$ \\
31079 & 2 & 2 & 4 & 3 & 1 & 2 & 5 & 690 & $D_{5}$ \\
31159 & 2 & 1 & 16 & 4 & 5 & 2 & 6 & 691 & $D_{3}$ \\
31247 & 2 & 6 & 4 & 3 & 1 & 2 & 3 & 692 & $D_{13}$ \\
31271 & 2 & 1 & 6 & 3 & 2 & 2 & 9 & 693 & $D_{3}$ \\
31321 & 2 & 3 & 8 & 4 & 2 & 2 & 7 & 693 & $D_{7}$ \\
31513 & 2 & 1 & 16 & 4 & 4 & 2 & 4 & 697 & $D_{3}$ \\
31543 & 2 & 2 & 4 & 3 & 1 & 2 & 5 & 697 & $D_{5}$ \\
31847 & 2 & 5 & 4 & 3 & 1 & 2 & 3 & 703 & $D_{11}$ \\
32009 & 2 & 1 & 12 & 4 & 3 & 2 & 9 & 706 & $D_{3}$ \\
32143 & 2 & 3 & 4 & 3 & 1 & 2 & 4 & 708 & $D_{7}$ \\
32183 & 2 & 3 & 4 & 3 & 1 & 2 & 3 & 708 & $D_{7}$ \\
32327 & 2 & 1 & 16 & 4 & 5 & 2 & 5 & 710 & $D_{3}$ \\
32327 & 2 & 3 & 4 & 3 & 1 & 2 & 3 & 710 & $D_{9}$ \\
32353 & 2 & 1 & 20 & 4 & 6 & 2 & 5 & 711 & $D_{3}$ \\
32401 & 2 & 1 & 20 & 4 & 6 & 2 & 6 & 712 & $D_{3}$ \\
32479 & 2 & 5 & 4 & 3 & 1 & 2 & 4 & 714 & $D_{11}$ \\
32647 & 2 & 3 & 4 & 3 & 1 & 2 & 4 & 719 & $D_{7}$ \\
32687 & 2 & 5 & 4 & 3 & 1 & 2 & 3 & 720 & $D_{11}$ \\
32719 & 2 & 2 & 4 & 3 & 1 & 2 & 14 & 721 & $D_{5}$ \\
32887 & 2 & 2 & 6 & 3 & 2 & 2 & 6 & 724 & $D_{5}$ \\
32983 & 2 & 1 & 4 & 3 & 1 & 2 & 4 & 725 & $D_{3}$ \\
33223 & 2 & 1 & 4 & 3 & 1 & 2 & 9 & 732 & $D_{3}$ \\
33343 & 2 & 1 & 4 & 3 & 1 & 2 & 5 & 733 & $D_{3}$ \\
33679 & 2 & 1 & 4 & 3 & 1 & 2 & 5 & 738 & $D_{3}$ \\
33767 & 2 & 3 & 4 & 3 & 1 & 2 & 3 & 739 & $D_{7}$ \\
34351 & 2 & 2 & 4 & 3 & 1 & 2 & 8 & 753 & $D_{5}$ \\
34471 & 2 & 2 & 12 & 4 & 3 & 2 & 6 & 756 & $D_{5}$ \\
34487 & 2 & 1 & 14 & 4 & 4 & 2 & 7 & 756 & $D_{3}$ \\
34591 & 2 & 2 & 4 & 3 & 1 & 2 & 5 & 757 & $D_{5}$ \\
34679 & 2 & 1 & 6 & 3 & 2 & 2 & 4 & 758 & $D_{3}$ \\
34679 & 2 & 3 & 4 & 3 & 1 & 2 & 4 & 758 & $D_{9}$ \\
34721 & 2 & 1 & 12 & 4 & 4 & 2 & 9 & 759 & $D_{3}$ \\
34847 & 2 & 1 & 16 & 4 & 5 & 2 & 6 & 762 & $D_{3}$ \\
35401 & 2 & 1 & 56 & 5 & 8 & 2 & 10 & 776 & $D_{3}$ \\
35591 & 2 & 9 & 4 & 3 & 1 & 2 & 7 & 779 & $D_{19}$ \\
35759 & 2 & 3 & 6 & 3 & 2 & 2 & 5 & 781 & $D_{7}$ \\
35839 & 2 & 2 & 4 & 3 & 1 & 2 & 5 & 781 & $D_{5}$ \\
35977 & 2 & 2 & 8 & 4 & 2 & 2 & 5 & 783 & $D_{5}$ \\
36191 & 2 & 1 & 46 & 5 & 7 & 2 & 11 & 786 & $D_{3}$ \\
36791 & 2 & 1 & 12 & 4 & 3 & 2 & 8 & 799 & $D_{3}$ \\
36871 & 2 & 2 & 6 & 3 & 2 & 2 & 5 & 801 & $D_{5}$ \\
37087 & 2 & 1 & 32 & 4 & 11 & 2 & 5 & 804 & $D_{3}$ \\
37199 & 2 & 1 & 18 & 4 & 4 & 2 & 8 & 806 & $D_{3}$ \\
37607 & 2 & 1 & 4 & 3 & 1 & 2 & 3 & 814 & $D_{3}$ \\
37831 & 2 & 2 & 4 & 3 & 1 & 2 & 5 & 820 & $D_{5}$ \\
37879 & 2 & 3 & 4 & 3 & 1 & 2 & 5 & 821 & $D_{7}$ \\
37993 & 2 & 2 & 8 & 4 & 2 & 2 & 5 & 824 & $D_{5}$ \\
38047 & 2 & 2 & 4 & 3 & 1 & 2 & 9 & 825 & $D_{5}$ \\
38167 & 2 & 1 & 14 & 4 & 4 & 2 & 10 & 829 & $D_{3}$ \\
38231 & 2 & 1 & 4 & 3 & 1 & 2 & 4 & 830 & $D_{3}$ \\
38287 & 2 & 1 & 4 & 3 & 1 & 2 & 5 & 832 & $D_{3}$ \\
38303 & 2 & 1 & 4 & 3 & 1 & 2 & 3 & 832 & $D_{3}$ \\
38593 & 2 & 1 & 20 & 4 & 6 & 2 & 9 & 836 & $D_{3}$ \\
38959 & 2 & 1 & 10 & 3 & 4 & 2 & 9 & 842 & $D_{3}$ \\
38977 & 2 & 1 & 20 & 4 & 6 & 2 & 13 & 842 & $D_{3}$ \\
39023 & 2 & 1 & 40 & 5 & 7 & 2 & 10 & 842 & $D_{3}$ \\
39199 & 2 & 2 & 4 & 3 & 1 & 2 & 7 & 844 & $D_{5}$ \\
39631 & 2 & 1 & 16 & 4 & 5 & 2 & 11 & 853 & $D_{3}$ \\
39679 & 2 & 1 & 6 & 3 & 2 & 2 & 12 & 854 & $D_{3}$ \\
39679 & 2 & 3 & 4 & 3 & 1 & 2 & 4 & 854 & $D_{9}$ \\
\end{longtable}

\subsection*{Characteristic $p=2$ and non-prime square-free levels, dihedral}

\begin{longtable}{||c|c|c|c|c|c|c|c|c|c||}
\hline
Level & Wt & ResD & Dim & EmbDim & NilO & GorDef & \#Ops & \#(p$<$HB) & Gp \\
\hline\endhead\hline\endfoot\hline\hline\endlastfoot
1055 & 2 & 1 & 24 & 4 & 9 & 2 & 5 & 47 & $D_{3}$ \\
1727 & 2 & 1 & 16 & 4 & 5 & 2 & 6 & 65 & $D_{3}$ \\
2071 {}\footnote{First found by W.\ Stein.} & 2 & 1 & 8 & 4 & 2 & 2 & 5 & 73 & $D_{3}$ \\
2631 & 2 & 1 & 40 & 4 & 17 & 2 & 6 & 106 & $D_{3}$ \\
2991 & 2 & 1 & 40 & 4 & 17 & 2 & 4 & 121 & $D_{3}$ \\
3095 & 2 & 1 & 40 & 4 & 17 & 2 & 4 & 114 & $D_{3}$ \\
3431 & 2 & 1 & 24 & 5 & 4 & 2 & 6 & 107 & $D_{3}$ \\
3471 & 2 & 1 & 16 & 5 & 3 & 2 & 5 & 146 & $D_{3}$ \\
3639 & 2 & 1 & 28 & 4 & 11 & 2 & 5 & 140 & $D_{3}$ \\
4031 & 2 & 1 & 16 & 4 & 4 & 2 & 6 & 125 & $D_{3}$ \\
4087 & 2 & 1 & 8 & 4 & 2 & 2 & 6 & 126 & $D_{3}$ \\
4119 & 2 & 1 & 12 & 4 & 3 & 2 & 4 & 156 & $D_{3}$ \\
4415 & 2 & 1 & 8 & 4 & 2 & 2 & 6 & 153 & $D_{3}$ \\
\end{longtable}

\subsection*{Characteristic $p=2$, icosahedral}

\begin{longtable}{||c|c|c|c|c|c|c|c|c|c||}
\hline
Level & Wt & ResD & Dim & EmbDim & NilO & GorDef & \#Ops & \#(p$<$HB) & Gp \\
\hline\endhead\hline\endfoot\hline\hline\endlastfoot
89491 & 2 & 2 & 12 & 4 & 3 & 2 & 4 & 1746 & $A_5$ \\
\end{longtable}

\subsection*{Characteristic $p=3$, prime levels, dihedral}

\begin{longtable}{||c|c|c|c|c|c|c|c|c|c||}
\hline
Level & Wt & ResD & Dim & EmbDim & NilO & GorDef & \#Ops & \#(p$<$HB) & Gp \\
\hline\endhead\hline\endfoot\hline\hline\endlastfoot
1031 & 3 & 2 & 4 & 3 & 1 & 2 & 4 & 55 & $D_{5}$ \\
1511 & 3 & 3 & 4 & 3 & 1 & 2 & 9 & 74 & $D_{7}$ \\
2087 & 3 & 2 & 4 & 3 & 1 & 2 & 3 & 98 & $D_{5}$ \\
4259 & 3 & 2 & 4 & 3 & 1 & 2 & 3 & 179 & $D_{5}$ \\
4799 & 3 & 3 & 22 & 4 & 9 & 2 & 9 & 196 & $D_{7}$ \\
5939 & 3 & 2 & 4 & 3 & 1 & 2 & 4 & 235 & $D_{5}$ \\
6899 & 3 & 2 & 4 & 3 & 1 & 2 & 3 & 269 & $D_{5}$ \\
6959 & 3 & 2 & 4 & 3 & 1 & 2 & 4 & 270 & $D_{5}$ \\
7523 & 3 & 2 & 4 & 3 & 1 & 2 & 4 & 289 & $D_{5}$ \\
7559 & 3 & 2 & 4 & 3 & 1 & 2 & 6 & 290 & $D_{5}$ \\
7583 & 3 & 3 & 20 & 3 & 9 & 2 & 6 & 290 & $D_{7}$ \\
8219 & 3 & 2 & 4 & 3 & 1 & 2 & 6 & 310 & $D_{5}$ \\
8447 & 3 & 5 & 20 & 3 & 9 & 2 & 3 & 318 & $D_{11}$ \\
8699 & 3 & 2 & 6 & 3 & 2 & 2 & 9 & 326 & $D_{5}$ \\
9431 & 3 & 3 & 4 & 3 & 1 & 2 & 4 & 350 & $D_{7}$ \\
9743 & 3 & 2 & 8 & 3 & 3 & 2 & 3 & 360 & $D_{5}$ \\
9887 & 3 & 2 & 8 & 3 & 3 & 2 & 3 & 365 & $D_{5}$ \\
10079 & 3 & 2 & 60 & 3 & 29 & 2 & 5 & 368 & $D_{5}$ \\
10247 & 3 & 2 & 10 & 4 & 3 & 2 & 5 & 375 & $D_{5}$ \\
10847 & 3 & 3 & 22 & 4 & 9 & 2 & 9 & 395 & $D_{7}$ \\
12011 & 3 & 2 & 4 & 3 & 1 & 2 & 3 & 431 & $D_{5}$ \\
12119 & 3 & 2 & 56 & 3 & 27 & 2 & 8 & 434 & $D_{5}$ \\
12263 & 3 & 2 & 8 & 3 & 3 & 2 & 3 & 438 & $D_{5}$ \\
12959 & 3 & 5 & 20 & 3 & 9 & 2 & 4 & 457 & $D_{11}$ \\
13907 & 3 & 2 & 22 & 4 & 9 & 2 & 8 & 487 & $D_{5}$ \\
14699 & 3 & 2 & 4 & 3 & 1 & 2 & 6 & 513 & $D_{5}$ \\
14783 & 3 & 3 & 20 & 3 & 9 & 2 & 3 & 515 & $D_{13}$ \\
14783 & 3 & 3 & 20 & 3 & 9 & 2 & 3 & 515 & $D_{13}$ \\
\end{longtable}

\subsection*{Characteristic $p=5$, prime levels, dihedral}

\begin{longtable}{||c|c|c|c|c|c|c|c|c|c||}
\hline
Level & Wt & ResD & Dim & EmbDim & NilO & GorDef & \#Ops & \#(p$<$HB) & Gp \\
\hline\endhead\hline\endfoot\hline\hline\endlastfoot
419 & 5 & 1 & 4 & 3 & 1 & 2 & 3 & 40 & $D_{3}$ \\
439 & 5 & 1 & 14 & 4 & 5 & 2 & 10 & 42 & $D_{3}$ \\
491 & 5 & 1 & 4 & 3 & 1 & 2 & 3 & 46 & $D_{3}$ \\
751 & 5 & 1 & 12 & 3 & 5 & 2 & 3 & 65 & $D_{3}$ \\
839 & 5 & 1 & 6 & 3 & 2 & 2 & 6 & 70 & $D_{3}$ \\
1231 & 5 & 1 & 4 & 3 & 1 & 2 & 3 & 97 & $D_{3}$ \\
2579 & 5 & 1 & 4 & 3 & 1 & 2 & 3 & 180 & $D_{3}$ \\
2699 & 5 & 1 & 14 & 4 & 5 & 2 & 8 & 188 & $D_{3}$ \\
3299 & 5 & 1 & 4 & 3 & 1 & 2 & 6 & 220 & $D_{3}$ \\
3359 & 5 & 1 & 4 & 3 & 1 & 2 & 6 & 222 & $D_{3}$ \\
4111 & 5 & 1 & 4 & 3 & 1 & 2 & 3 & 267 & $D_{3}$ \\
4219 & 5 & 1 & 20 & 3 & 6 & 2 & 5 & 274 & $D_{3}$ \\
4931 & 5 & 3 & 12 & 3 & 5 & 2 & 3 & 310 & $D_{7}$ \\
5011 & 5 & 1 & 4 & 3 & 1 & 2 & 5 & 316 & $D_{3}$ \\
5639 & 5 & 1 & 6 & 3 & 2 & 2 & 5 & 348 & $D_{3}$ \\
5939 & 5 & 3 & 12 & 3 & 5 & 2 & 5 & 366 & $D_{7}$ \\
6079 & 5 & 1 & 6 & 3 & 2 & 2 & 5 & 370 & $D_{3}$ \\
6271 & 5 & 1 & 4 & 3 & 1 & 2 & 3 & 379 & $D_{3}$ \\
6571 & 5 & 1 & 12 & 3 & 5 & 2 & 5 & 399 & $D_{3}$ \\
6691 & 5 & 1 & 4 & 3 & 1 & 2 & 5 & 405 & $D_{3}$ \\
6779 & 5 & 1 & 6 & 3 & 2 & 2 & 6 & 410 & $D_{3}$ \\
7459 & 5 & 1 & 12 & 3 & 5 & 2 & 7 & 443 & $D_{3}$ \\
7759 & 5 & 3 & 4 & 3 & 1 & 2 & 3 & 457 & $D_{7}$ \\
8779 & 5 & 1 & 12 & 3 & 5 & 2 & 12 & 511 & $D_{3}$ \\
8819 & 5 & 3 & 4 & 3 & 1 & 2 & 3 & 513 & $D_{7}$ \\
9011 & 5 & 1 & 4 & 3 & 1 & 2 & 6 & 522 & $D_{3}$ \\
\end{longtable}

\subsection*{Characteristic $p=7$, prime levels, dihedral}

\begin{longtable}{||c|c|c|c|c|c|c|c|c|c||}
\hline
Level & Wt & ResD & Dim & EmbDim & NilO & GorDef & \#Ops & \#(p$<$HB) & Gp \\
\hline\endhead\hline\endfoot\hline\hline\endlastfoot
199 & 7 & 1 & 4 & 3 & 1 & 2 & 3 & 30 & $D_{3}$ \\
839 & 7 & 1 & 4 & 3 & 1 & 2 & 6 & 93 & $D_{3}$ \\
1259 & 7 & 1 & 4 & 3 & 1 & 2 & 3 & 130 & $D_{3}$ \\
1291 & 7 & 1 & 4 & 3 & 1 & 2 & 4 & 133 & $D_{3}$ \\
1319 & 7 & 1 & 4 & 3 & 1 & 2 & 6 & 136 & $D_{3}$ \\
1399 & 7 & 1 & 4 & 3 & 1 & 2 & 3 & 141 & $D_{3}$ \\
1559 & 7 & 1 & 4 & 3 & 1 & 2 & 7 & 155 & $D_{3}$ \\
1567 & 7 & 1 & 8 & 3 & 4 & 2 & 3 & 156 & $D_{3}$ \\
1823 & 7 & 1 & 6 & 3 & 2 & 2 & 4 & 179 & $D_{3}$ \\
1823 & 7 & 3 & 4 & 3 & 1 & 2 & 4 & 179 & $D_{9}$ \\
\end{longtable}

\subsection*{Characteristic $p=11$, prime levels, dihedral}

\begin{longtable}{||c|c|c|c|c|c|c|c|c|c||}
\hline
Level & Wt & ResD & Dim & EmbDim & NilO & GorDef & \#Ops & \#(p$<$HB) & Gp \\
\hline\endhead\hline\endfoot\hline\hline\endlastfoot
431 & 11 & 3 & 4 & 3 & 1 & 2 & 5 & 77 & $D_{7}$ \\
563 & 11 & 1 & 4 & 3 & 1 & 2 & 3 & 97 & $D_{3}$ \\
1187 & 11 & 1 & 6 & 3 & 2 & 2 & 6 & 181 & $D_{3}$ \\
1223 & 11 & 3 & 4 & 3 & 1 & 2 & 4 & 187 & $D_{7}$ \\
1231 & 11 & 1 & 4 & 3 & 1 & 2 & 3 & 189 & $D_{3}$ \\
1231 & 11 & 3 & 4 & 3 & 1 & 2 & 3 & 189 & $D_{9}$ \\
1327 & 11 & 1 & 4 & 3 & 1 & 2 & 3 & 199 & $D_{5}$ \\
1327 & 11 & 1 & 4 & 3 & 1 & 2 & 4 & 199 & $D_{5}$ \\
1583 & 11 & 1 & 24 & 3 & 11 & 2 & 4 & 230 & $D_{3}$ \\
1619 & 11 & 1 & 4 & 3 & 1 & 2 & 3 & 235 & $D_{3}$ \\
1823 & 11 & 1 & 4 & 3 & 1 & 2 & 4 & 263 & $D_{3}$ \\
2243 & 11 & 1 & 4 & 3 & 1 & 2 & 3 & 310 & $D_{3}$ \\
2351 & 11 & 1 & 4 & 3 & 1 & 2 & 6 & 325 & $D_{3}$ \\
2351 & 11 & 3 & 4 & 3 & 1 & 2 & 6 & 325 & $D_{9}$ \\
2503 & 11 & 1 & 4 & 3 & 1 & 2 & 3 & 341 & $D_{3}$ \\
2591 & 11 & 1 & 4 & 3 & 1 & 2 & 5 & 351 & $D_{3}$ \\
2647 & 11 & 1 & 4 & 3 & 1 & 2 & 3 & 360 & $D_{3}$ \\
2767 & 11 & 1 & 4 & 3 & 1 & 2 & 3 & 370 & $D_{3}$ \\
2791 & 11 & 1 & 4 & 3 & 1 & 2 & 3 & 375 & $D_{3}$ \\
3011 & 11 & 1 & 4 & 3 & 1 & 2 & 5 & 402 & $D_{3}$ \\
3119 & 11 & 1 & 4 & 3 & 1 & 2 & 5 & 415 & $D_{3}$ \\
3299 & 11 & 1 & 4 & 3 & 1 & 2 & 4 & 434 & $D_{3}$ \\
3299 & 11 & 3 & 4 & 3 & 1 & 2 & 4 & 434 & $D_{9}$ \\
3571 & 11 & 1 & 4 & 3 & 1 & 2 & 4 & 462 & $D_{3}$ \\
\end{longtable}

\subsection*{Characteristic $p=13$, prime levels, dihedral}

\begin{longtable}{||c|c|c|c|c|c|c|c|c|c||}
\hline
Level & Wt & ResD & Dim & EmbDim & NilO & GorDef & \#Ops & \#(p$<$HB) & Gp \\
\hline\endhead\hline\endfoot\hline\hline\endlastfoot
367 & 13 & 1 & 4 & 3 & 1 & 2 & 3 & 78 & $D_{3}$ \\
439 & 13 & 2 & 4 & 3 & 1 & 2 & 4 & 91 & $D_{5}$ \\
563 & 13 & 1 & 4 & 3 & 1 & 2 & 3 & 111 & $D_{3}$ \\
971 & 13 & 2 & 4 & 3 & 1 & 2 & 4 & 177 & $D_{5}$ \\
1223 & 13 & 2 & 4 & 3 & 1 & 2 & 4 & 216 & $D_{5}$ \\
1427 & 13 & 1 & 4 & 3 & 1 & 2 & 5 & 243 & $D_{3}$ \\
1439 & 13 & 1 & 28 & 3 & 13 & 2 & 5 & 246 & $D_{3}$ \\
1823 & 13 & 1 & 4 & 3 & 1 & 2 & 4 & 298 & $D_{3}$ \\
\end{longtable}

\subsection*{Characteristic $p=17$, prime levels, dihedral}

\begin{longtable}{||c|c|c|c|c|c|c|c|c|c||}
\hline
Level & Wt & ResD & Dim & EmbDim & NilO & GorDef & \#Ops & \#(p$<$HB) & Gp \\
\hline\endhead\hline\endfoot\hline\hline\endlastfoot
59 & 17 & 1 & 6 & 3 & 2 & 2 & 3 & 23 & $D_{3}$ \\
239 & 17 & 2 & 4 & 3 & 1 & 2 & 5 & 68 & $D_{5}$ \\
1327 & 17 & 1 & 4 & 3 & 1 & 2 & 3 & 289 & $D_{3}$ \\
1427 & 17 & 2 & 4 & 3 & 1 & 2 & 3 & 306 & $D_{5}$ \\
1951 & 17 & 1 & 4 & 3 & 1 & 2 & 4 & 402 & $D_{3}$ \\
2503 & 17 & 1 & 4 & 3 & 1 & 2 & 3 & 497 & $D_{3}$ \\
2687 & 17 & 1 & 36 & 3 & 17 & 2 & 4 & 529 & $D_{3}$ \\
\end{longtable}

\subsection*{Characteristic $p=19$, prime levels, dihedral}

\begin{longtable}{||c|c|c|c|c|c|c|c|c|c||}
\hline
Level & Wt & ResD & Dim & EmbDim & NilO & GorDef & \#Ops & \#(p$<$HB) & Gp \\
\hline\endhead\hline\endfoot\hline\hline\endlastfoot
439 & 19 & 1 & 6 & 3 & 2 & 2 & 3 & 125 & $D_{3}$ \\
751 & 19 & 1 & 4 & 3 & 1 & 2 & 5 & 195 & $D_{5}$ \\
751 & 19 & 1 & 4 & 3 & 1 & 2 & 5 & 195 & $D_{5}$ \\
1427 & 19 & 1 & 6 & 3 & 2 & 2 & 3 & 335 & $D_{3}$ \\
\end{longtable}

\subsection*{Characteristic $p=23$, prime levels, dihedral}

\begin{longtable}{||c|c|c|c|c|c|c|c|c|c||}
\hline
Level & Wt & ResD & Dim & EmbDim & NilO & GorDef & \#Ops & \#(p$<$HB) & Gp \\
\hline\endhead\hline\endfoot\hline\hline\endlastfoot
83 & 23 & 1 & 4 & 3 & 1 & 2 & 3 & 37 & $D_{3}$ \\
503 & 23 & 3 & 4 & 3 & 1 & 2 & 4 & 162 & $D_{7}$ \\
971 & 23 & 2 & 4 & 3 & 1 & 2 & 4 & 284 & $D_{5}$ \\
1259 & 23 & 1 & 4 & 3 & 1 & 2 & 3 & 358 & $D_{3}$ \\
\end{longtable}

\subsection*{Characteristic $p=29$, prime levels, dihedral}

\begin{longtable}{||c|c|c|c|c|c|c|c|c|c||}
\hline
Level & Wt & ResD & Dim & EmbDim & NilO & GorDef & \#Ops & \#(p$<$HB) & Gp \\
\hline\endhead\hline\endfoot\hline\hline\endlastfoot
107 & 29 & 1 & 4 & 3 & 1 & 2 & 3 & 55 & $D_{3}$ \\
199 & 29 & 1 & 4 & 3 & 1 & 2 & 3 & 92 & $D_{3}$ \\
\end{longtable}

\subsection*{Characteristic $p=31$, prime levels, dihedral}

\begin{longtable}{||c|c|c|c|c|c|c|c|c|c||}
\hline
Level & Wt & ResD & Dim & EmbDim & NilO & GorDef & \#Ops & \#(p$<$HB) & Gp \\
\hline\endhead\hline\endfoot\hline\hline\endlastfoot
367 & 31 & 1 & 4 & 3 & 1 & 2 & 3 & 161 & $D_{3}$ \\
743 & 31 & 3 & 4 & 3 & 1 & 2 & 4 & 293 & $D_{7}$ \\
\end{longtable}

\subsection*{Characteristic $p=37$, prime levels, dihedral}

\begin{longtable}{||c|c|c|c|c|c|c|c|c|c||}
\hline
Level & Wt & ResD & Dim & EmbDim & NilO & GorDef & \#Ops & \#(p$<$HB) & Gp \\
\hline\endhead\hline\endfoot\hline\hline\endlastfoot
139 & 37 & 1 & 4 & 3 & 1 & 2 & 4 & 83 & $D_{3}$ \\
\end{longtable}

\subsection*{Characteristic $p=41$, prime levels, dihedral}

\begin{longtable}{||c|c|c|c|c|c|c|c|c|c||}
\hline
Level & Wt & ResD & Dim & EmbDim & NilO & GorDef & \#Ops & \#(p$<$HB) & Gp \\
\hline\endhead\hline\endfoot\hline\hline\endlastfoot
83 & 41 & 1 & 4 & 3 & 1 & 2 & 3 & 61 & $D_{3}$ \\
139 & 41 & 1 & 4 & 3 & 1 & 2 & 4 & 92 & $D_{3}$ \\
431 & 41 & 1 & 4 & 3 & 1 & 2 & 5 & 233 & $D_{3}$ \\
\end{longtable}

\subsection*{Characteristic $p=43$, prime levels, dihedral}

\begin{longtable}{||c|c|c|c|c|c|c|c|c|c||}
\hline
Level & Wt & ResD & Dim & EmbDim & NilO & GorDef & \#Ops & \#(p$<$HB) & Gp \\
\hline
419 & 43 & 1 & 4 & 3 & 1 & 2 & 3 & 239 & $D_{3}$ \\
\hline
\end{longtable}

\subsection*{Characteristic $p=47$, prime levels, dihedral}

\begin{longtable}{||c|c|c|c|c|c|c|c|c|c||}
\hline
Level & Wt & ResD & Dim & EmbDim & NilO & GorDef & \#Ops & \#(p$<$HB) & Gp \\
\hline\endhead\hline\endfoot\hline\hline\endlastfoot
31 & 47 & 1 & 4 & 3 & 1 & 2 & 3 & 30 & $D_{3}$ \\
107 & 47 & 1 & 4 & 3 & 1 & 2 & 3 & 82 & $D_{3}$ \\
139 & 47 & 1 & 4 & 3 & 1 & 2 & 4 & 101 & $D_{3}$ \\
179 & 47 & 2 & 4 & 3 & 1 & 2 & 3 & 126 & $D_{5}$ \\
\end{longtable}

\subsection*{Characteristic $p=53$, prime levels, dihedral}

\begin{longtable}{||c|c|c|c|c|c|c|c|c|c||}
\hline
Level & Wt & ResD & Dim & EmbDim & NilO & GorDef & \#Ops & \#(p$<$HB) & Gp \\
\hline\endhead\hline\endfoot\hline\hline\endlastfoot
131 & 53 & 2 & 4 & 3 & 1 & 2 & 3 & 106 & $D_{5}$ \\
211 & 53 & 1 & 4 & 3 & 1 & 2 & 4 & 159 & $D_{3}$ \\
\end{longtable}

\subsection*{Characteristic $p=59$, prime levels, dihedral}

\begin{longtable}{||c|c|c|c|c|c|c|c|c|c||}
\hline
Level & Wt & ResD & Dim & EmbDim & NilO & GorDef & \#Ops & \#(p$<$HB) & Gp \\
\hline\endhead\hline\endfoot\hline\hline\endlastfoot
23 {}\footnote{First found by K.\ Buzzard, unpublished.} & 59 & 1 & 4 & 3 & 1 & 2 & 4 & 30 & $D_{3}$ \\
211 & 59 & 1 & 4 & 3 & 1 & 2 & 4 & 175 & $D_{3}$ \\
227 & 59 & 1 & 4 & 3 & 1 & 2 & 3 & 187 & $D_{5}$ \\
227 & 59 & 1 & 4 & 3 & 1 & 2 & 3 & 187 & $D_{5}$ \\
367 & 59 & 1 & 4 & 3 & 1 & 2 & 3 & 279 & $D_{3}$ \\
\end{longtable}

\subsection*{Characteristic $p=61$, prime levels, dihedral}

\begin{longtable}{||c|c|c|c|c|c|c|c|c|c||}
\hline
Level & Wt & ResD & Dim & EmbDim & NilO & GorDef & \#Ops & \#(p$<$HB) & Gp \\
\hline\endhead\hline\endfoot\hline\hline\endlastfoot
239 & 61 & 1 & 4 & 3 & 1 & 2 & 4 & 199 & $D_{5}$ \\
239 & 61 & 1 & 4 & 3 & 1 & 2 & 4 & 199 & $D_{5}$ \\
431 & 61 & 1 & 4 & 3 & 1 & 2 & 5 & 327 & $D_{3}$ \\
\end{longtable}

\subsection*{Characteristic $p=67$, prime levels, dihedral}

\begin{longtable}{||c|c|c|c|c|c|c|c|c|c||}
\hline
Level & Wt & ResD & Dim & EmbDim & NilO & GorDef & \#Ops & \#(p$<$HB) & Gp \\
\hline\endhead\hline\endfoot\hline\hline\endlastfoot
31 & 67 & 1 & 4 & 3 & 1 & 2 & 6 & 41 & $D_{3}$ \\
239 & 67 & 2 & 4 & 3 & 1 & 2 & 5 & 217 & $D_{5}$ \\
\end{longtable}

\subsection*{Characteristic $p=71$, prime levels, dihedral}

\begin{longtable}{||c|c|c|c|c|c|c|c|c|c||}
\hline
Level & Wt & ResD & Dim & EmbDim & NilO & GorDef & \#Ops & \#(p$<$HB) & Gp \\
\hline\endhead\hline\endfoot\hline\hline\endlastfoot
59 & 71 & 1 & 4 & 3 & 1 & 2 & 3 & 71 & $D_{3}$ \\
239 & 71 & 1 & 4 & 3 & 1 & 2 & 5 & 223 & $D_{3}$ \\
283 & 71 & 1 & 4 & 3 & 1 & 2 & 5 & 263 & $D_{3}$ \\
\end{longtable}

\subsection*{Characteristic $p=73$, prime levels, dihedral}

\begin{longtable}{||c|c|c|c|c|c|c|c|c|c||}
\hline
Level & Wt & ResD & Dim & EmbDim & NilO & GorDef & \#Ops & \#(p$<$HB) & Gp \\
\hline\endhead\hline\endfoot\hline\hline\endlastfoot
211 & 73 & 1 & 4 & 3 & 1 & 2 & 4 & 209 & $D_{3}$ \\
283 & 73 & 1 & 4 & 3 & 1 & 2 & 5 & 269 & $D_{3}$ \\
\end{longtable}

\subsection*{Characteristic $p=79$, prime levels, dihedral}

\begin{longtable}{||c|c|c|c|c|c|c|c|c|c||}
\hline
Level & Wt & ResD & Dim & EmbDim & NilO & GorDef & \#Ops & \#(p$<$HB) & Gp \\
\hline\endhead\hline\endfoot\hline\hline\endlastfoot
307 & 79 & 1 & 4 & 3 & 1 & 2 & 5 & 307 & $D_{3}$ \\
\end{longtable}

\subsection*{Characteristic $p=83$, prime levels, dihedral}

\begin{longtable}{||c|c|c|c|c|c|c|c|c|c||}
\hline
Level & Wt & ResD & Dim & EmbDim & NilO & GorDef & \#Ops & \#(p$<$HB) & Gp \\
\hline\endhead\hline\endfoot\hline\hline\endlastfoot
47 & 83 & 2 & 4 & 3 & 1 & 2 & 4 & 67 & $D_{5}$ \\
79 & 83 & 2 & 4 & 3 & 1 & 2 & 3 & 101 & $D_{5}$ \\
107 & 83 & 1 & 6 & 3 & 2 & 2 & 3 & 132 & $D_{3}$ \\
211 & 83 & 1 & 4 & 3 & 1 & 2 & 4 & 232 & $D_{3}$ \\
251 & 83 & 1 & 4 & 3 & 1 & 2 & 3 & 271 & $D_{7}$ \\
251 & 83 & 1 & 4 & 3 & 1 & 2 & 3 & 271 & $D_{7}$ \\
251 & 83 & 1 & 4 & 3 & 1 & 2 & 3 & 271 & $D_{7}$ \\
\end{longtable}

\subsection*{Characteristic $p=89$, prime levels, dihedral}

\begin{longtable}{||c|c|c|c|c|c|c|c|c|c||}
\hline
Level & Wt & ResD & Dim & EmbDim & NilO & GorDef & \#Ops & \#(p$<$HB) & Gp \\
\hline\endhead\hline\endfoot\hline\hline\endlastfoot
131 & 89 & 1 & 4 & 3 & 1 & 2 & 3 & 165 & $D_{5}$ \\
131 & 89 & 1 & 4 & 3 & 1 & 2 & 3 & 165 & $D_{5}$ \\
\end{longtable}

\subsection*{Characteristic $p=97$, prime levels, dihedral}

\begin{longtable}{||c|c|c|c|c|c|c|c|c|c||}
\hline
Level & Wt & ResD & Dim & EmbDim & NilO & GorDef & \#Ops & \#(p$<$HB) & Gp \\
\hline\endhead\hline\endfoot\hline\hline\endlastfoot
307 & 97 & 1 & 4 & 3 & 1 & 2 & 5 & 367 & $D_{3}$ \\
\end{longtable}

\end{document}